\documentclass[12pt,amscd]{amsart}
\usepackage{amsmath,amsxtra,amssymb,latexsym, amscd,amsthm}
\usepackage{indentfirst}
\usepackage[mathscr]{eucal}
\usepackage[dvips]{color}
\usepackage{graphicx}
\usepackage{graphics,epsfig}
\usepackage{tikz,pgf,pgfplots}
\usepackage{mathrsfs}
\usepackage{graphicx,color}
\usetikzlibrary{arrows}
\usepackage{multirow}
\usepackage{verbatim}

\usepackage[pagebackref=true]{hyperref}

\newtheorem{Theorem}{Theorem}[section]
\newtheorem{Definition}[Theorem]{Definition}
\newtheorem{Lemma}[Theorem]{Lemma}
\newtheorem{Corollary}[Theorem]{Corollary}
\newtheorem{Proposition}[Theorem]{Proposition}
\theoremstyle{definition}
\newtheorem{Example}[Theorem]{Example}
\newtheorem*{Remark}{Remark}

\def\NZQ{\Bbb}

\def\ZZ{{\NZQ Z}}

\def\C{\mathfrak C}
\def\F {\mathcal F}

\def\girth{\operatorname{girth}}
\def\h {\widetilde{H}}

\def\mm{{\frak m}}
\def\a{{\bold a}}
\def\b{{\bold b}}

\def\0{{\bold 0}}
\def\x{{\bold x}}
\def\al{{\alpha}}
\def\D{{\Delta}}

\DeclareMathOperator{\CS}{CS}
\DeclareMathOperator{\Gin}{Gin}

\begin{document}
\title {Extremal Betti numbers of certain two-dimensional monomial ideals}

\author{Nguy\^en Quang L\^oc}
\address{Department of Mathematics, Hanoi National University of Education, 136 Xuan Thuy, Hanoi, Vietnam}
\email{nqloc@hnue.edu.vn}

\author{Nguy\^en C\^ong Minh}
\address{Faculty of Mathematics and Informatics, Hanoi University of Science and Technology, 1 Dai Co Viet, Hanoi, Vietnam}
\email{minh.nguyencong@hust.edu.vn}
\email{ngcminh@gmail.com}

\author{Phan Thi Thuy}
\address{Department of Mathematics, Hanoi National University of Education, 136 Xuan Thuy, Hanoi, Vietnam}
\email{phanthuy@hnue.edu.vn}  
\date{}          
\maketitle

 \let\thefootnote\relax
  
\begin{abstract}
In this paper, we shall provide explicit formulas for the extremal Betti numbers of $R/I$, where $I$ is the defining ideal of certain weighted hyperplanes in $\Bbb{P}^{n-1}$ and $R$ is the polynomial ring in $n$ indeterminates over a field. As a consequence, we completely classify such ideals which are pseudo-Gorenstein as in sense of V. Ene, J. Herzog, T. Hibi and S. S. Madani \cite{EHHM} .
\end{abstract}
\section{Introduction} 
Let $R = K[x_1,\ldots, x_n]$ be a polynomial ring over a field $K$ with the maximal homogeneous ideal $\mm = (x_1,\ldots,x_n)$. Let $M$ be a finitely generated graded $R$-module.  Let
$$
0 \longrightarrow \bigoplus_{j \geq 0}R(-j)^{\beta_{p,j}} \longrightarrow \bigoplus_{j \geq 0}R(-j)^{\beta_{p-1,j}} \longrightarrow \cdots \longrightarrow \bigoplus_{j \geq 0}R(-j)^{\beta_{0,j}} \longrightarrow M \longrightarrow 0
$$
be the minimal graded free resolution of $M$, where $p$ is its projective dimension. Here, $R(-j)$ is the free $R$-module of rank 1 generated in degree $j$ and $\beta_{i,j} = \beta_{i,j}(M)$ is the $(i, j)$-th graded Betti number of $M$.  The graded Betti numbers can be arranged in a so-called Betti diagram, which is a table whose element in the $i$-th column and the $j$-th row is $\beta_{i, i+j}$. Such an $(i, j)$-position is called a corner, and the corresponding graded Betti number $\beta_{i, i+j}$ called extremal, if 
$$
\beta_{i, i +j} \neq 0 \quad\text{ and }\quad \beta_{r, r +s} = 0 \text{ for all pairs $(r, s) \neq (i, j)$ with $r \geq i$ and $s \geq j$}.
$$

The extremal Betti numbers of the graded $R$-module $M$, introduced by D. Bayer, H. Charalambous and S. Popescu in \cite[Definition 4.3.13]{BCP}, can be seen as a refinement of the Castelnuovo-Mumford regularity and of the projective dimension of  $M$. Projective dimension measures the column index of the easternmost extremal Betti number, whereas regularity measures the row index of the southernmost extremal Betti number. Moreover, the extremal Betti numbers of $M$, as well as their positions, are preserved by passing from $M$ to $\Gin(M)$, the generic initial module with respect to the graded reverse lexicographic order in the case $K$ is a field of characteristic 0. 


Both position and value of any extremal Betti number can be read off from the local cohomology modules with support in $\mm$ by a result of P. Schenzel \cite{Sch} (see also \cite{Chd}) which relates the extremal Betti numbers to the extremal local cohomologies. Similar to the Betti diagram, the local cohomology diagram of $M$ is a table whose element in the $i$-th column and the $j$-th row is $\dim_K H^i_{\mm}(M)_{j-i}$, where $H^{i}_{\mm}(M)$ is the $i$-th local cohomology module with support in $\mm$ of $M$. An $(i, j)$-position of this diagram is called a corner, and the corresponding local cohomology is called extremal, if
$$
H^i_{\mm}(M)_{j - i} \neq 0 \text{ and } H^u_{\mm}(M)_{v - u} = 0 \text{ for all $(u, v) \neq (i, j)$ with $u \leq i$ and $v \geq j$}.
$$

\begin{Theorem}[\cite{Sch}] \label{corner}
A pair $(i, j)$ is a corner of the Betti diagram if and only if $(n-i, j)$ is a corner of the local cohomology diagram of $M$. In this case, we have 
$$
\beta_{i, i+j}(M) = \dim_K H^{n-i}_{\mm}(M)_{i + j - n}.
$$
\end{Theorem}

Recently monomial ideals have attracted a lot of attention, both for their own sake
and (for instance via generic initial ideals) as a tool for studying more general ideals, and problems in combinatorics and geometry. Many questions remain open, however. One important question is finding the Betti numbers of a monomial ideal. To get an explicit formula for Betti numbers is difficult, even extremal Betti numbers (see \cite{BCP,EHHM,G, HR, HSV,LMT}).

Denote by $\C_n$ the set of the ideals $I$ in $R$ of the form
$$ 
I = \bigcap_{1\leq i<j\leq n}P_{i,j}^{w_{i,j}},
$$
where $P_{i,j}$ are ideals of $R$ generated by variables $\{x_1,x_2,...,x_n\}\backslash \{x_i,x_j\}$ for $1\leq i<j\leq n$ and $w_{i,j}$ is a non-negative integer. 
Let $T$ be the union of $\binom{n}{2}$ hyperplanes $L_{i,j}$ defined by $x_i=x_j=0$ in $\Bbb{P}^{n-1}$. The number $w_{ij}$ is usually called the weight of the hyperplane $L_{i,j}$. So $I$ can be viewed as the defining ideal of $T$ with respect to the weights $w_{i,j}$. 

The starting point of this work has been Migliore-Nagel's paper \cite{MNa},  where they considered $n=4$ and the question of whether or not these ideals were Cohen-Macaulay. After that, some other algebraic properties of these ideals have been intensively studied, see  \cite{GH, MN, MT, Th}. In \cite{MN, MT}, the Cohen-Macaulayness of $R/I$ and its $a_i$-invariants were investigated when $w_{i,j}$ only takes either $\alpha$ or $\beta$ and $n\ge 5$, where the $a_i$-invariant of the $\ZZ$-graded module $M$, $a_i(M)$, is defined by the top non-vanishing degree of $H^{i}_{\mm}(M)$ if it is non-zero, otherwise $a_i(M)=- \infty$. It is noted that if $\alpha=\beta$ or $\alpha>0=\beta$, such an ideal $I$ can be viewed as a symbolic power of a squarefree monomial ideal. In this case, the extremal Betti numbers of $R/I$ are recently determined in \cite{LMT}. It is worth mentioning that in \cite{Ma}, the author obtained the Betti table of the symbolic power of a star configuration of hypersurfaces. It means that the Betti table of the defining ideal of $T$ with respect to the same weight $w_{i,j}$ can be explicitly computed.

Let $\C_n(\alpha,\beta)$ consist of $I\in\C_n$ with $w_{i,j} \in \{\alpha, \beta\}$ for all $i, j$, where $\alpha>\beta>0$ are given. Our aim in this article is to determine the extremal Betti numbers of $R/I$ for  $I$ in $\C_n(\alpha,\beta)$. 

Since $\dim(R/I) = 2$ and $H^0_{\mm}(R/I) = 0$, the local cohomology diagram of $R/I$ has at most two corners, which are $(1, a_1(R/I) + 1)$ and $(2, a_2(R/I) + 2)$. Consequently, by Theorem~\ref{corner}, the possible extremal Betti numbers of $R/I$ are
$$
\beta_{n-1, n + a_1(R/I)} \quad \text{ and } \quad \beta_{n-2, n + a_2(R/I)}.
$$
In particular, if $R/I$ is Cohen-Macaulay, then $a_1(R/I) = -\infty$ and $\beta_{n-2, n + a_2(R/I)}$ is the unique extremal Betti number. On the other hand, if $R/I$ is not Cohen-Macaulay, then $\beta_{n-1, n + a_1(R/I)}$ is the unique extremal Betti number if and only if $(1, a_1(R/I) + 1)$ is the unique corner of the local cohomology diagram. This happens only when $a_1(R/I) + 1 \geq a_2(R/I) + 2$. The extremal Betti numbers are determined via dimensions of the extremal local cohomologies by Theorem~\ref{corner}, which will be in turn calculated in terms of reduced homology groups by Takayama's formula (Lemma~\ref{L1}).

We now present the layout of the article. In Section 2 we recall some basic notions and results which are needed, including Takayama's formula mentioned above. Our main results, Theorems~\ref{extremal1}--\ref{extremal5}, are gathered in Section 3. In Section 4, we compute the values of $\beta_{n-2,n + a_2(R/I)}$ when it is extremal in Proposition~\ref{3}. Section 5 deals with the values of $\beta_{n-1,n + a_1(R/I)}$ when it is extremal and $\alpha \geq \beta + 3$, with the main result being Proposition~\ref{thm2}. Together, Propositions~\ref{3} and \ref{thm2} imply Theorem~\ref{extremal2}, which is the case when $\alpha \geq \beta + 3$. When $\alpha = \beta + 1$, the extremal Betti numbers of $R/I$ are computed in Theorem~\ref{extremal1}, and when $\alpha = \beta + 2$ they are determined in Theorems \ref{extremal3}, \ref{extremal4} and \ref{extremal5}. In Section 6, we classify the rings $R/I$ that are pseudo-Gorenstein as a consequence of our theorems and provide numerous examples illustrating our results.

\section{Preliminaries}
In this section, we recall some definitions and properties concerning simplicial complexes and graph theory that will be used later. The interested reader is referred to \cite{S}, \cite{D} for more details. 

\subsection{Simplicial complexes and Stanley-Reisner correspondence}  A simplicial complex $\D$ on the finite set $[n]=\{1,2,\ldots,n\}$ is a collection of subsets of $[n]$ called faces, closed under taking subsets; that is, if $F\in \D$ is a face and $G\subseteq F$, then $G\in \D$. If a face $F\in \D$ has the cardinality $\vert F\vert = i +1$, then $i$ is called the dimension of $F$. The dimension of $\D$ is the maximum of the dimensions of its faces, or it is $-\infty$ if $\D$ is the void complex, which has no faces. It is clear that $\D$ can be uniquely determined by the set of its maximal elements under inclusion, called by facets, which is denoted $\F(\D)$. The Stanley-Reisner ideal of $\D$ is the following ideal of $R = K[x_1, \ldots, x_n]$:
$$
I_{\D} = (x_{i_1}\cdots x_{i_s} \mid \{i_1, \ldots, i_s\} \not \in \D) = \bigcap_{F \in \F(\D)} P_F,
$$
where $P_F$ is the prime ideal of $R$ generated by all variables $x_i$ with $i \not \in F$. The ideal $I_{\D}$ is a squarefree monomial ideal. Conversely, if $I$ is a squarefree monomial ideal of $R$, then it is the Stanley-Reisner ideal associated to the simplicial complex 
$$
\D(I) = \{\{i_1, \ldots, i_s\} \subseteq [n] \mid x_{i_1}\cdots x_{i_s} \not \in I\}.
$$

\subsection{Simplicial homology} Let $\Delta$ be a simplicial complex. An oriented $q$-simplex of $\Delta$ is a face $F \in \Delta$, $|F| = q+1$, with an ordering of the vertices, with the rule that two orderings define the same orientation if and only if they differ by an even permutation. Let $C_q(\Delta)$ be the $K$-vector space with basis consisting of the oriented $q$-simplices of $\Delta$. We define the homomorphisms $\partial_q:C_q(\Delta) \to C_{q-1}(\Delta)$ for $q \ge 1$ by defining them on the basis elements by 
$$\partial_q [v_0, \ldots,v_q] = \sum_{i=0}^q (-1)^i [v_0, \ldots, \hat v_i, \ldots, v_q],$$
where $\hat v_i$ denotes that $v_i$ is missing. It is easily verified that $\partial_q \partial_{q+1} = 0$. The chain complex $C(\Delta) = \{C_q(\Delta),\partial_q\}$ is called the oriented chain complex of $\Delta$. Let $C_{-1}(\Delta)$ be the $K$-vector space with basis $\{\emptyset\}$, and define an augmentation $\epsilon:C_0(\Delta) \to C_{-1}(\Delta)$ by $\epsilon(x) = \emptyset$ for every $x \in [n]$. The augmented chain complex $(C(\Delta),\epsilon)$ is called the augmented oriented chain complex of $\Delta$.

\begin{Definition} The $q$-th reduced homology group of $\Delta$ with coefficients in $K$, denoted $\h_q(\Delta; K)$, is defined to be the $q$-th homology group of the augmented oriented chain complex of $\Delta$ over $K$.
\end{Definition}

A simplicial complex $\D$ is called {\it acyclic} if $\h_i(\Delta;K) = 0$ for all $i$.
It is noted that if $\Delta$ is the empty complex (i.e., $\Delta=\{\emptyset\}$), then $\h_i(\Delta;K) \neq 0$ if and only if $i = -1$. Also, if $\Delta$ is a cone over some $x \in [n]$ or $\Delta$ is the void complex  (i.e., $\Delta=\emptyset$), then it is acyclic.

\subsection{Degree complexes} Let $J$ be a monomial ideal in $R$. We set $\D(J) = \D(\sqrt{J})$. Inspired by a result of Hochster in the squarefree case \cite[Theorem 4.1]{S}, Takayama found the following combinatorial formula for $\dim_K H_\mm^i (R/J)_\a$ \cite[Proposition 1]{T}. For a subset $F$ of $[n]$, let $R_F=R[x_i^{-1}\mid i\in F]$. For $\a=(a_1,\ldots,a_n)\in\mathbb{Z}^n$, the co-support of $\a$ is defined to be the set $\CS_{\a} = \{i \in [n] \mid a_i <0\}$. Let
$$
\D_\a(J) =\{F \subseteq [n]\setminus \CS_{\a} \mid \x^\a = x_1^{a_1}\cdots x_n^{a_n}\not\in JR_{F\cup \CS_{\a}} \}.
$$

\begin{Lemma}[Takayama's formula] \label{L1}
For all $\a\in \mathbb{Z}^n$ and $i\geq 0$, we have
\begin{equation*}
\dim_K H^i_\mm(R/J)_\a=
\begin{cases}
\dim_K \h_{i-\vert \CS_{\a}\vert -1}(\D_\a(J);K) &\text{\quad  if $\CS_{\a}\in \D(J)$,}\\
0 &\text{\quad  otherwise.}
\end{cases}
\end{equation*}
\end{Lemma}

For $\a = (a_1,\ldots,a_n)\in\mathbb{Z}^n$, we put $\vert \a\vert =a_1+\cdots+a_n$ and $\sigma_{i,j}(\a)=\vert \a\vert - a_i-a_j$, where $1\leq i<j\leq n$.

\begin{Lemma}[\cite{MN}, Lemmas 2.2, 2.3] \label{L2} Let $I\in\C_n(\al,\beta)$. For any $\a\in \mathbb{N}^n$ we have 
\begin{enumerate}
\item[(i)] $\dim \D_\a(I)\leq 1$.
\item[(ii)] $\{i,j\}\in\D_\a(I)$ if and only if $\sigma_{i,j}(\a)<w_{i,j}$ for any $1\leq i<j\leq n$.
\item[(iii)] If $\D_\a(I)\neq \{\emptyset\}$, then $\dim \D_\a(I)=1$ and $\D_\a(I)$ is pure (i.e., $\Delta _\a(I)$ can be reviewed as a graph for all $\a\in\mathbb{N}^n$).
\end{enumerate}
\end{Lemma}

\begin{Lemma}[\cite{MN}, Proposition 2.4] \label{L3} Let $I\in\C_n(\al,\beta)$ and $\a\in \mathbb{Z}^n$. Then $H^1_\mm(R/I)_\a \neq 0$ if and only if $\a \in \mathbb{N}^n$ (i.e., $\CS_{\a} = \emptyset$) and $\D_\a(I)$ is a disconnected graph.
\end{Lemma}

\subsection{Graphs} A simplicial complex of dimension one can be considered as a simple graph, i.e., a graph without multiple edges or loops. Let $G$ be a (simple) graph with the vertex set $V(G) \subseteq [n]$ and the edge set $E(G)$. For simplicity, an edge ${\{i, j\} \in E(G)}$ will also be written as $ij \in G$, while the notation $ij \not \in G$ means $i \neq j$ and $\{i, j\}$ is not an edge. For each $i \in [n]$, the neighborhood of $i$ in $G$ is $
N(i)=\{j\in V(G)\mid ij \in G\},
$
and let
$
N[i] = N(i) \cup \{i\}.
$

The degree of $i$, denoted $\deg(i)$, is the cardinality of the set $N(i)$. A cycle $C_t = \{12, 23, \ldots, t1\}$ of length $t \geq 3$, or a $t$-cycle, is written as $(1, 2, \ldots, t)$. The girth of $G$, denoted $\girth(G)$, is the smallest length of cycles of $G$; if $G$ contains no cycles, we set $\girth(G) = \infty$.

The next lemma is well-known (e.g., \cite[Exercises 5.2.7]{V}) and will be used implicitly throughout the article.

\begin{Lemma} \label{Euler}
Let $G$ be a simple graph considered as a simplicial complex of dimension one, with $c(G)$ connected components. Then 
$$
\dim_K \widetilde H_0(G;K) = c(G) -1 \quad \text{ and } \quad \dim_K \widetilde H_1(G;K) = |E(G)| + c(G) - |V(G)|.
$$
\end{Lemma}


\section{The main results}

Let $I$ be an ideal in $\C_n(\alpha,\beta)$ for $n \geq 5$. We define
$$
G = \{\{i,j\}\mid w_{i,j}=\alpha, 1\leq i<j\leq n\}.
$$
We note that $G$ is a non-void simple graph with the vertex set $V(G)\subseteq [n]$ having no isolated vertices and the edge set $E(G)$. 
Thus, $i \in V(G)$ if and only if there exists $j \in [n]$ such that $ij$ is an edge of $G$ (then $j \in V(G)$ as well). For a subset $W$ of $[n]$, let $$G[W]=\{ij\in G \mid i,j \in W\}$$ denote the induced subgraph of $G$ by $W$. The graph $G$ plays an important role since the invariants of $R/I$ will be computed in terms of $G$. 

\medskip

Following \cite{LMT}, a pair of disjoint edges $(ij, pq)$ of $G$ is said to be

- {\it disconnected}, if the induced subgraph $G[i,j,p,q]$ is $\{ij, pq\}$;

- {\it type-1}, if the induced subgraph $G[i,j,p,q]$ is a path of length 3;

- {\it type-2}, if the induced subgraph $G[i,j,p,q]$ is a path of length 4 which is not a 4-cycle. Namely, it is a triangle with a whisker.

We denote by $p_0(G)$ the number of disconnected pairs of edges, $p_1(G)$ the number of type-1 pairs of edges, and $p_2(G)$ the number of type-2 pairs of edges of $G$, respectively. 

We define an invariant of the graph $G$, denoted $a(G)$, by
\begin{align}
a(G) =\underset{ij\not\in G}\sum|N(i)\cap N(j)|.
\end{align}
Our first main result determines the extremal Betti numbers of $R/I$ when $\alpha = \beta + 1$.

\begin{Theorem} \label{extremal1}
Let $\alpha =\beta +1$. The extremal Betti numbers of $R/I$ are:
\begin{center}
\begin{tabular}{ |p{2.5cm}|p{3cm}|p{3cm}|p{6.5cm}| }
\hline
Corners of Betti diagram & \multicolumn{2}{c|}{Structure of $G$} & Extremal Betti numbers of $R/I$ \\ \hline 
 &  & \vspace*{0.0001cm}$\beta = 3,$ $s(G)\neq 0$  & $\beta_{n-1,n+5}(R/I)=  s(G)$; \\
&& & $\beta_{n-2, n + 9}(R/I) =c_3(G)$.\\
\cline{3-4} & $\girth(G) = 3$ &$\beta$ is odd,   & $\beta_{n-1,n+\alpha+\beta -2}(R/I)=  (n-4)p_0(G);$\\
&& $\beta \neq 1,3,$&\\
&&$p_0(G)\neq 0$ & $\beta_{n-2, n + 3\alpha-3} (R/I)=c_3(G)$.\\
\cline{3-4} & &\vspace*{0.0001cm}$\beta$ is even;  &$\beta_{n-1,n+\alpha+\beta -1}(R/I)=  p_0(G);$\\
&& $p_0(G)~\neq~ 0$& $\beta_{n-2, n + 3\alpha-3}(R/I) =c_3(G)$.\\
\cline{2-4} & &  \vspace*{0.02cm} $\beta = 3;$ $s(G)\neq 0$  &  $\beta_{n-1, n +5}(R/I) =  s(G)$;
\\
 \vspace*{0.1cm} $(n-1, a_1 + 1);$ && &$\beta_{n-2, n + 8}(R/I) =
\underset{1\leq i \leq n}\sum \binom{\deg(i)}{2}.$
\\ 
\cline{3-4}   & $ \girth(G) \neq 3,$  $G$   & $\beta$ is odd, &  $\beta_{n-1,n+\alpha+\beta -2}(R/I)=  (n-4)p_0(G);$\\
$(n-2, a_2 + 2)$&contains a vertex & $\beta \neq 1,3,$&\\
 &of degree $\geq 2$& $p_0(G)\neq 0$ & $\beta_{n-2, n + 2\alpha+\beta-3}(R/I) =\underset{1\leq i \leq n}\sum \binom{\deg(i)}{2}.$\\

\cline{3-4} &&\vspace*{0.01cm}$\beta$ is even; &$\beta_{n-1,n+\alpha+\beta -1}(R/I)=  p_0(G);$\\
&&$p_0(G)~\neq~ 0$  & $\beta_{n-2, n + 2\alpha+\beta-3}(R/I) =\underset{1\leq i \leq n}\sum \binom{\deg(i)}{2}.$\\ 

\cline{2-4}&  & \vspace*{0.01cm}$\beta = 3$  & $\beta_{n-1,n+5}(R/I)=  s(G)$; \\
&& & $\beta_{n-2, n + 7}(R/I) =(n-2)|E(G)|.$\\

\cline{3-4}& $G$ consists of disjoint edges; & \vspace*{0.0001cm}$\beta$ is odd,  & $\beta_{n-1,n+\alpha+\beta -2}(R/I)=  (n-4)p_0(G);$\\
 & $|E(G)|\neq 1$&$\beta ~\neq~ 1, 3$& $\beta_{n-2, n + \alpha+2\beta-3}(R/I) =(n-2)|E(G)|.$\\
\cline{3-4} & &\vspace*{0.01cm}$\beta$ is even &$\beta_{n-1,n+\alpha+\beta -1}(R/I)=  p_0(G);$\\
&&& $\beta_{n-2, n + \alpha+2\beta-3}(R/I) =(n-2)|E(G)|.$\\
 \hline

 &\vspace*{0.01cm}$\girth(G) = 3$&  & \vspace*{0.01cm}
 $\beta_{n-2, n + 3\alpha-3}(R/I) = c_3(G).$
\\
\cline{2-2}\cline{4-4}\vspace*{0.6cm}$(n-2, a_2 + 2)$& \vspace*{0.1cm}$ \girth(G) \neq 3;$ $G$ contains a vertex of degree $\geq 2$ &   \vspace*{0.01cm}$\beta = 1;$ or \qquad\qquad\qquad $\beta =3, s(G)= 0;$ or $\beta ~\neq~ 1,3,$  $p_0(G)=0$ & \vspace*{0.01cm}$\beta_{n-2, n + 2\alpha+\beta-3}(R/I) = \qquad\quad\quad\quad$
\hspace*{1.5cm}$\begin{cases}
a(G)    &\text{ if } \beta =1, \\
\underset{1\leq i \leq n}\sum \binom{\deg(i)}{2} &\text{ if } \beta \neq 1.
\end{cases}
$
\\ 
\cline{2-2}\cline{4-4}& $G$ consists of disjoint edges & &\vspace*{0.001cm} $\beta_{n-2, n + \alpha+2\beta-3}(R/I) =(n-2)|E(G)|.$\\
\hline
\end{tabular}
\end{center}
\vspace{0.2cm}
where $s(G) = |\{W\subseteq [n] \mid |W|=5 \text{ and } G[W] \text{ is a disconnected graph}\}|$ and 
$c_3(G)~$~is~the number of 3-cycles in $G.$ 

\end{Theorem}
\newpage

In the case $\alpha \geq \beta + 2$, we consider two conditions on the graph $G$, namely
\begin{enumerate}
\item[$(G_1)$] There exists a pair of disjoint edges of $G$ which is not contained in any cycle of length 4.
\item[$(G_2)$] There exists $ij \not \in G$ such that $\{i, j\} \cup G$ is a disconnected graph (where $\{i, j\}$ is an edge of that graph).
\end{enumerate}

If $R/I$ is not Cohen-Macaulay, then $G$ must satisfy either the condition $(G_1)$ or the condition $(G_2)$ by \cite[Theorem 4.5]{MN}. Observe that if $G$ does not satisfy the condition $(G_1)$, then $G$ is a connected graph.

We define another invariant of $G$, which is

\begin{align}
b(G) &= {\underset{pq\in G}\sum} \, |V(G) \backslash (N(p)\cup N(q))|.
\end{align}


\begin{Remark} It is easy to verify the following observations:
\begin{enumerate}
\item [(i)] The graph $G$ satisfies the condition $(G_1)$ if and only if $p_0(G)+p_1(G)+p_2(G)$ is a positive number. It is also equivalent to the condition $b(G) > 0$.

\item [(ii)] If $G$ is a connected graph, then the condition $(G_2)$ is equivalent to
\begin{align*}
|\{ij \not \in G \mid ij\cup G \text{ is disconnected}\}| 
= \binom{n-|V(G)|}{2} > 0.
\end{align*}
\end{enumerate}
\end{Remark}
\noindent For an integer $m$, we set  
$$
\overline{m} = \max\{m, 0\},
$$
and denote 
\begin{align*}
g_1(\alpha,\beta) & = \biggl[(\alpha - \beta)^2 + (\alpha - \beta - 1)^2 - \overline{\alpha - 2\beta + 1}^2\biggr]p_0(G) + \\
& \dfrac{1}{2}\biggl[{(\alpha - \beta - 1)(2\alpha - 2\beta - 1)} - {\overline{\alpha - 2\beta + 1}^2}\biggr] p_1(G) + \\
& \dfrac{1}{4}\biggl[2{(\alpha - \beta - 1)^2} - {\overline{\alpha - 2\beta + 1}^2}\biggr] p_2(G) + \overline{\alpha - 2\beta}\, b(G); \qquad\qquad
\end{align*}
\begin{align*}
g_2(\alpha,\beta) &= \biggl[(\alpha - \beta)^2 + (\alpha - \beta - 1)^2 - \overline{\alpha - 2\beta + 1}^2\biggr]p_0(G) + \\
& \dfrac{1}{2}\biggl[(\alpha - \beta)^2 + (\alpha - \beta - 1)(\alpha - \beta - 2) - \overline{\alpha - 2\beta + 1}^2\biggr] p_1(G) + \\
& \dfrac{1}{4}\biggl[(\alpha - \beta)^2 + (\alpha - \beta - 2)^2 - \overline{\alpha - 2\beta + 1}^2\biggr] p_2(G) + \overline{\alpha - 2\beta}\, b(G); 
\end{align*}
\begin{align*}
g_3(\alpha,\beta) &= \biggl[2(\alpha - \beta)(\alpha - \beta - 1) - \overline{\alpha - 2\beta}\, (\alpha - 2\beta +2)\biggr]p_0(G) + \\
& \dfrac{1}{2}\biggl[{(\alpha - \beta)(2\alpha - 2\beta - 3)} - \overline{\alpha - 2\beta}\, (\alpha - 2\beta + 2)\biggr] p_1(G) + \\
& \dfrac{1}{4}\biggl[2{(\alpha - \beta)(\alpha - \beta - 2)} - \overline{\alpha - 2\beta}\, (\alpha - 2\beta + 2)\biggr] p_2(G) + \overline{\alpha - 2\beta}\, b(G); 
\end{align*}
\begin{align*}
g_4(\alpha,\beta) &= \biggl[2(\alpha - \beta)(\alpha - \beta - 1) - \overline{\alpha - 2\beta}\, (\alpha - 2\beta +2)\biggr]p_0(G) + \\
& \dfrac{1}{2}\biggl[{(\alpha - \beta - 1)(2\alpha - 2\beta - 1)} - \overline{\alpha - 2\beta}\, (\alpha - 2\beta + 2)\biggr] p_1(G) + \\
& \dfrac{1}{4}\biggl[2(\alpha - \beta - 1)^2 - \overline{\alpha - 2\beta}\, (\alpha - 2\beta + 2)\biggr] p_2(G) + \overline{\alpha - 2\beta}\, b(G);\qquad\qquad
\end{align*}
\begin{align*}
 h_1(\alpha, \beta) = \dfrac{1}{2}\biggl[{(\alpha - \beta - 1)^2}- {\overline{\alpha - 2\beta + 1}^2}+ 2\, \overline{\alpha - 2\beta}\, |V(G)|\biggr]\binom{n-|V(G)|}{2};\qquad \qquad 
\end{align*}
\begin{align*}
h_2(\alpha, \beta) = \dfrac{1}{4}\biggl[(\alpha - \beta)^2 + (\alpha - \beta - 2)^2 - 2\, \overline{\alpha - 2\beta + 1}^2 + 2\, \overline{\alpha - 2\beta}\, |V(G)|\biggr] \binom{n-|V(G)|}{2};  
\end{align*}
\begin{align*}
 h_3(\alpha, \beta) = \dfrac{1}{2}\biggl[{(\alpha - \beta)(\alpha - \beta - 2)}- \overline{\alpha - 2\beta}\, (\alpha - 2\beta + 2)+ 2\, \overline{\alpha - 2\beta}\, |V(G)|\biggr]\binom{n-|V(G)|}{2}; 
\end{align*}
\begin{align*}
h_4(\alpha, \beta) = \dfrac{1}{2}\biggl[{(\alpha - \beta - 1)^2}- \overline{\alpha - 2\beta}\, (\alpha - 2\beta + 2)+ 2\, \overline{\alpha - 2\beta}\, |V(G)|\biggr]\binom{n-|V(G)|}{2}. 
\end{align*}
We set
\begin{equation*}
B_1(\alpha,\beta)=
\begin{cases}
g_1(\alpha,\beta) & \text{ if $\alpha$ is odd and $\beta$ is even;}\\
g_2(\alpha,\beta) & \text{ if $\alpha$ is odd and $\beta$ is odd;}\\
g_3(\alpha,\beta) &  \text{ if $\alpha$ is even and $\beta$ is even;}\\
g_4(\alpha,\beta) & \text{ if $\alpha$ is even and $\beta$ is odd;}\\
\end{cases} 
\end{equation*}

and

\begin{equation*}
B_2(\alpha,\beta)=
\begin{cases}
h_1(\alpha,\beta) & \text{ if $\alpha$ is odd and $\beta$ is even;}\\
h_2(\alpha,\beta) & \text{ if $\alpha$ is odd, $\beta$ is odd and $\beta > 1$;}\\
h_3(\alpha,\beta) & \text{ if $\alpha$ is even and $\beta$ is even;}\\
h_4(\alpha,\beta) & \text{ if $\alpha$ is even, $\beta$ is odd and $\beta > 1$;}\\
(\alpha - 2)\binom{n-|V(G)|}{2} & \text{ if $\beta=1$}.
\end{cases}
\end{equation*}

\medskip

Here, $B_1(\alpha, \beta)$ refers to the formula of $\beta_{n-1, n + a_1(R/I)}$ when $G$ satisfies $(G_1)$, while $B_2(\alpha, \beta)$ refers to $\beta_{n-1, n + a_1(R/I)}$ when $G$ does not satisfy $(G_1)$ (so it must satisfies $(G_2)$).
We are now able to describe the extremal Betti numbers of $R/I$ when $\alpha \geq \beta + 3$.
\begin{Theorem} \label{extremal2}
Let $\alpha \geq \beta + 3$. The extremal Betti numbers of $R/I$ are:
\begin{center}
\begin{tabular}{ |p{2.5cm}|p{2.5cm}|p{3cm}|p{6.5cm}| }
\hline
Corners of Betti diagram & \multicolumn{2}{c|}{Structure of $G$} & Extremal Betti numbers of $R/I$ \\ \hline 
 &  & \vspace*{0.0001cm}$G$ satisfies $(G_1)$  & $\beta_{n-1,n+2\alpha -2}(R/I)=  B_1(\alpha,\beta)$; \\
&& & $\beta_{n-2, n + 3\alpha-3}(R/I) =c_3(G)$.\\
\cline{3-4} &$\girth(G) = 3$ &\vspace*{0.0001cm}$G$ satisfies $(G_2)$; &$\beta_{n-1,n+\alpha+\beta -2}(R/I) =  B_2(\alpha,\beta)$; \\
& & not $(G_1)$ &   $\beta_{n-2, n + 3\alpha-3}(R/I) =c_3(G)$.\\
\cline{2-4} & &&$\beta_{n-1, n + 2\alpha-2}(R/I)= B_1(\alpha,\beta);$
\\
 \vspace*{0.5cm} \vspace*{0.5cm} $(n-1, a_1 + 1);$ & $ \girth(G) \neq 3,$ $G$ contains a vertex of degree $\geq 2$  &\vspace*{0.35cm} $G$ satisfies $(G_1)$&\vspace{0.1cm}$\beta_{n-2, n + 2\alpha+\beta-3}(R/I) = \qquad\qquad\qquad$
\hspace*{0.8cm}$\begin{cases}
a(G)  + (\alpha - 2)c_4(G)  \text{ if } \beta =1,\\
\underset{1\leq i \leq n}\sum \binom{\deg(i)}{2}\quad\qquad \text{ } \text{ }\text{ if } \beta \neq 1.
\end{cases}
$
\\ 
\cline{3-4} $(n-2, a_2 + 2)$ &   & & $\beta_{n-1, n + \alpha+\beta-2}(R/I) = B_2(\alpha,\beta);$\\

 &&$G$ satisfies $(G_2)$; not $(G_1)$& \vspace{0.1cm}$\beta_{n-2, n + 2\alpha+\beta-3}(R/I) = \qquad\qquad\qquad$
\hspace*{0.8cm}$\begin{cases}
a(G)  + (\alpha - 2)c_4(G)  \text{ if } \beta =1, \\
\underset{1\leq i \leq n}\sum \binom{\deg(i)}{2}\quad\qquad \text{ } \text{ }\text{ if } \beta \neq 1.
\end{cases}
$
\\
\cline{2-4}\cline{3-3}&  \vspace{0.01cm}$G$ consists of  &  \vspace{0.2cm}$|E(G)| = 1$ & $\beta_{n-1, n + \alpha+\beta-2}(R/I) = B_2(\alpha,\beta);$\\
  &disjoint edges&&$\beta_{n-2, n + \alpha+2\beta-3}(R/I) =n-2.$\\

\cline{3-4}&&\vspace{0.1cm}$|E(G)| \neq 1,$ & $\beta_{n-1, n + 2\alpha-2}(R/I) =B_1(\alpha,\beta);$ \\
&&$\alpha<2\beta$  &  \vspace{0.01cm}$\beta_{n-2, n + \alpha+2\beta-3}(R/I) =(n-2)|E(G)|.$\\
 \hline
 \vspace*{0.01cm}$(n-1, a_1 + 1)$ &  $G$ consists of disjoint edges & $|E(G)| \neq 1,$ \text{ }\qquad  $\alpha\geq 2\beta$&  \vspace*{0.01cm}$\beta_{n-1, n + 2\alpha-2}(R/I) = B_1(\alpha,\beta). $
\\
\hline
& $ \girth(G) = 3$ &  &  $\beta_{n-2, n + 3\alpha-3}(R/I) =c_3(G)$.\\ 

\cline{2-2}\cline{4-4} \vspace*{0.1cm}$(n-2, a_2 + 2)$ &  \vspace*{0.4cm}$ \girth(G) \neq 3$ &$G$ does not satisfy both $(G_1)$ and $(G_2)$ &
 $\beta_{n-2, n + 2\alpha+\beta-3}(R/I) =\qquad\qquad\qquad$
\hspace*{0.8cm}$\begin{cases}
a(G)  + (\alpha - 2)c_4(G)  \text{ if } \beta =1, \\
\underset{1\leq i \leq n}\sum \binom{\deg(i)}{2}\quad\qquad\text{ } \text{ } \text{ if } \beta \neq 1.
\end{cases}
$
\\
\hline
\end{tabular}
\end{center}
\vspace{0.2cm}
where  $c_4(G)$ is the number of 4-cycles in $G$.
\end{Theorem}

When $\alpha = \beta + 2$, there are three different cases to be considered: $\beta$ is odd; $\beta$ is even and $\beta \neq 2$; and $(\alpha, \beta) = (4, 2)$.

\begin{Theorem} \label{extremal3}
Let $\alpha = \beta + 2$ and $\beta$ is odd. The extremal Betti numbers of $R/I$ are:
\begin{center}
\begin{tabular}{ |p{2.5cm}|p{2.5cm}|p{3cm}|p{6.5cm}| }
\hline
Corners of Betti diagram & \multicolumn{2}{c|}{Structure of $G$} & Extremal Betti numbers of $R/I$ \\ \hline 
 &  & \vspace*{0.7cm}$G$ satisfies $(G_1)$ & 
$\beta_{n-1,n+\alpha+\beta}(R/I)=\qquad\qquad\quad\qquad$
$\begin{cases}
p_0(G)+b(G) \qquad\qquad\quad\text{ if } \beta=1, \\
5p_0(G)+2p_1(G)+p_2(G) \text{ if } \beta \geq 3;
\end{cases}$\\
& &&\vspace*{0.001cm} $\beta_{n-2, n + 3\alpha-3}(R/I) =c_3(G)$.\\
\cline{3-4} &$\girth(G) = 3$  & \vspace*{0.2cm}$G$ satisfies $(G_2)$, not $(G_1)$ &$\beta_{n-1,n+\alpha+\beta-2}(R/I)= \qquad\qquad\quad\quad$
\hspace*{1cm}$\begin{cases}
\binom{n-|V(G)|}{2} &\text{ if } \beta=1, \\
|E(G)|\binom{n-|V(G)|}{2} &\text{ if } \beta \geq 3;
\end{cases}$\\
&  & & \vspace*{0.001cm}  $\beta_{n-2, n + 3\alpha-3}(R/I) =c_3(G)$.\\

\cline{2-4}&  &\vspace*{1cm} $G$ satisfies $(G_1)$ & 
$\beta_{n-1,n+\alpha+\beta}(R/I)=\qquad\qquad\quad\qquad$
\hspace*{1cm}$\begin{cases}
p_0(G)+b(G) & \text{ if } \beta=1, \\
5p_0(G)+2p_1(G) & \text{ if } \beta \geq 3;
\end{cases}$\\
 \vspace*{0.5cm}$(n-1,a_1+1);$& &&\vspace*{0.001cm} $\beta_{n-2, n + 2\alpha+\beta-3}(R/I) =\qquad\qquad\quad\quad$
\hspace*{1cm}$\begin{cases}
a(G)  + c_4(G) & \quad\text{ if } \beta =1, \\
\underset{1\leq i \leq n}\sum \binom{\deg(i)}{2}& \quad\text{ if } \beta \geq 3.
\end{cases}$\\

\cline{3-4}$(n-2,a_2+2)$ &$\girth(G) \neq 3,$ $G$ contains a vertex of degree $\geq 2$  & \vspace*{0.5cm}$G$ satisfies $(G_2)$, not $(G_1)$ &$\beta_{n-1,n+\alpha+\beta-2}(R/I)= \qquad\qquad\quad\quad$
\hspace*{1cm}$\begin{cases}
\binom{n-|V(G)|}{2} & \text{ if } \beta=1, \\
|E(G)|\binom{n-|V(G)|}{2}& \text{ if } \beta \geq 3;
\end{cases}$\\
 &  & &   $\beta_{n-2, n + 2\alpha+\beta-3}(R/I) =\qquad\qquad\quad\quad$
\hspace*{1cm}$\begin{cases}
a(G)  + c_4(G) & \quad\text{ if } \beta =1, \\
\underset{1\leq i \leq n}\sum \binom{\deg(i)}{2}& \quad\text{ if } \beta \geq 3.
\end{cases}$\\

\cline{2-4} &   & $|E(G)|\neq 1,$  & 
$\beta_{n-1,n+\alpha+\beta}(R/I)=5p_0(G).$\\

&$G$ consists of & $\beta\geq 3$& $\beta_{n-2, n + \alpha+2\beta-3}(R/I) = (n-2)|E(G)|.$\\
\cline{3-4} & disjoint edges & $|E(G)|= 1$ & $\beta_{n-1,n+\alpha+\beta-2}(R/I)= \binom{n-2}{2};$ \\
&  & & \vspace*{0.0001cm}  $\beta_{n-2, n + \alpha+2\beta-3}(R/I) = n-2.$\\ \hline
$(n-1,a_1+1)$ & G consists of disjoint edges & $|E(G)|\neq 1,$ \text{ }$\beta =1$ & \vspace*{0.1cm}$\beta_{n-1,n+\alpha+\beta}(R/I)=p_0(G)+b(G). $\\
\hline

 & $\girth(G) = 3$ & &
$\beta_{n-2, n + 3\alpha-3}(R/I) =c_3(G)$.\\
\cline{2-2}\cline{4-4}\vspace*{0.2cm}$(n-2,a_2+2)$  &\vspace*{0.2cm} $\girth(G) \neq 3$  & $G$ does not satisfy both $(G_1)$ and $(G_2)$& $\beta_{n-2,n+2\alpha+\beta-3}(R/I)= \qquad\qquad\quad\quad$
\hspace*{1cm}$\begin{cases}
a(G)  + c_4(G) & \quad\text{ if } \beta =1, \\
\underset{1\leq i \leq n}\sum \binom{\deg(i)}{2}&\quad\text{ if } \beta \geq 3.
\end{cases}$\\

\hline
\end{tabular}
\end{center}
\end{Theorem}

In the case $\alpha = \beta + 2$ and $\beta$ is even, we need to consider the following conditions on $G$: 
\begin{enumerate}
\item[$(G_3)$] For any pair of disjoint edges $(ij, pq)$ of $G$ which is not contained in any cycle of length $4$, the induced subgraph $G[i,j,p,q]$ has form $(F_1)$. 
\item[$(G_4)$] If $G[i,j,p,q]$ has form $(F_1)$, then for every vertex $t\in [n]\setminus\{i,j,p,q\}$, the induced subgraph $G[i,j,p,q,t]$ has a subgraph that has form $(F_2).$
\item[$(G_5)$] $\{i,j\}\cup\{i,t\}\cup\{t,j\}\cup G$ is a connected graph for every $ij\not \in G$ and $ t\in[n]\backslash \{i,j\}$.\\
\\
\vspace{0.3cm}
\\
\\

\begin{center}
\setlength{\unitlength}{1.5cm}
\begin{picture}(5,1.5)\thicklines
\put(0,0){\line(0,1){1}}
\put(0,0){\line(1,0){2}}
\put(2,0){\line(0,0){1}}
\put(0,0){\line(2,1){2}}
\put(0,0){\circle*{0.1}}
\put(0,1){\circle*{0.1}}
\put(2,1){\circle*{0.1}}
\put(2,0){\circle*{0.1}}
\put(0.4,-1){Figure $(F_1)$}
\end{picture}
\begin{picture}(3,0)\thicklines
\put(0,0){\line(0,1){1}}
\put(0,0){\line(2,0){2}}
\put(2,0){\line(0,1){1}}
\put(0,0){\line(2,1){2}}
\put(2,0){\line(-1,2){1}}
\put(1,2){\circle*{0.1}}
\put(2,1){\line(-1,1){1}}
\put(0,1){\line(1,1){1}}
\put(0,0){\circle*{0.1}}
\put(0,1){\circle*{0.1}}
\put(2,1){\circle*{0.1}}
\put(2,0){\circle*{0.1}}
\put(1,2.1){t}
\put(0.4,-1){Figure $(F_2)$}
\end{picture}
\end{center}
\end{enumerate}

\vspace{2cm}

By \cite[Theorem 4.5 and Theorem 4.6]{MN}, $R/I$ is Cohen-Macaulay if and only if $G$  satisfies the conditions $(G_3), (G_4)$ and $(G_5)$ in the case $\beta=2,$ and $G$ does not satisfy both of the conditions $(G_1), (G_2)$  in other cases.

\medskip
Observe that when $G$ satisfies $(G_1)$, then it satisfies $(G_3)$ if and only if $$p_2(G) > 0 \text{ and } p_0(G) = p_1(G) = 0.$$

Let $f(G)$ denote the number of subgraphs $G[i, j, p, q, t]$ of $G$ with a subgraph of the form $(F_2)$, where $(ij, pq)$ is a type-2 pair of disjoint edges. We note that when $G$ satisfies both $(G_1)$ and $(G_3)$, the number $(n-4)p_2(G) - f(G)$ appeared in the below table (when $\beta = 4$) also equals 
$$
\underset{\underset{i\in V(G)\backslash (N(p)\cup N(q))}{pq\in G,} }\sum \left((n-4)\deg(i) - \binom{\deg(i)}{2}\right).
$$

Finally, we determine the extremal Betti numbers of $R/I$ in the case $\alpha = \beta + 2$ and $\beta$ is even by the following two theorems.
\newpage
\begin{Theorem} \label{extremal4}
Let $\alpha = \beta + 2$ and assume that $\beta$ is even, $(\alpha,\beta)\neq (4,2)$. The extremal Betti numbers of $R/I$ are:
\begin{center}
\begin{tabular}{ |p{2.5cm}|p{2.5cm}|p{3cm}|p{6.5cm}| }
\hline
Corners of Betti diagram & \multicolumn{2}{c|}{Structure of $G$} & Extremal Betti numbers of $R/I$ \\ \hline 
 &  & $G$ satisfies $(G_1),$ not $(G_3)$ & $\beta_{n-1,n+\alpha+\beta}(R/I)= 4p_0(G)+p_1(G);$\\
& && $\beta_{n-2, n + 3\alpha-3}(R/I) =c_3(G)$.\\

\cline{3-4} & \vspace*{1cm}$\girth(G) = 3$  & \vspace*{0.1cm}$G$ satisfies both  of $(G_1),$  $(G_3)$ & $\beta_{n-1,n+\alpha+\beta-1}(R/I)= \qquad\qquad\quad\quad$
$\begin{cases}
(n-4)p_2(G)-f(G) & \text{ if } \beta=4, \\
(n-4)p_2(G) & \text{ if } \beta \geq 6;
\end{cases}$\\
&  & & \vspace{0.01cm}  $\beta_{n-2, n + 3\alpha-3} (R/I)=c_3(G)$.\\

\cline{3-4} & & \vspace*{0.8cm}$G$ satisfies $(G_2),$ not $(G_1)$ &$\beta_{n-1,n+\alpha+\beta-3}(R/I)= \qquad\qquad\quad\quad$
$\begin{cases}
[(n-4)|E(G)|-\underset{1\leq i \leq n}\sum \binom{\deg(i)}{2}\\
\quad+c_3(G)]\binom{n-|V(G)|}{2} \qquad \text{ if } \beta = 4, \\
(n-4)|E(G)|\binom{n-|V(G)|}{2} \; \, \text{ if } \beta \geq 6;
\end{cases}$\\
&  & &  \vspace*{0.0001cm} $\beta_{n-2, n + 3\alpha-3}(R/I) =c_3(G)$.\\

\cline{2-4}$(n-1,a_1+1;)$ & &  $G$ satisfies $(G_1)$  & $\beta_{n-1,n+\alpha+\beta}(R/I)= 4p_0(G)+p_1(G);$\\
$(n-2,a_2+2)$ & $\girth(G)\neq 3,$  && \vspace{0.1cm}$\beta_{n-2, n + 2\alpha+\beta-3}(R/I) =\underset{1\leq i \leq n}\sum \binom{\deg(i)}{2}.$\\

\cline{3-4} &$G$ contains a vertex of degree $\geq 2$ &\vspace*{0.8cm}$G$ satisfies $(G_2),$ not $(G_1)$ &$\beta_{n-1,n+\alpha+\beta-3}(R/I)= \qquad\qquad\quad\quad$
$\begin{cases}
[(n-4)|E(G)|-\underset{1\leq i \leq n}\sum \binom{\deg(i)}{2}].\\
\hspace*{2.5cm}\binom{n-|V(G)|}{2} \text{ if } \beta=4, \\
(n-4)|E(G)|\binom{n-|V(G)|}{2} \text{ if } \beta \geq 6;
\end{cases}$\\
&  & & \vspace*{0.0001cm}  $\beta_{n-2, n + 2\alpha+\beta-3}(R/I) =\underset{1\leq i \leq n}\sum \binom{\deg(i)}{2}.$\\
\cline{2-4} &  & $|E(G)|\neq 1$ & 
$\beta_{n-1,n+\alpha+\beta}(R/I)= 4p_0(G);$\\
&\vspace*{0.1cm} $G$ consists of && \vspace*{0.0001cm}$\beta_{n-2, n + \alpha+2\beta-3}(R/I) = (n-2)|E(G)|.$\\
\cline{3-4} & disjoint edges & $|E(G)|= 1$ &$\beta_{n-1,n+\alpha+\beta-3}(R/I)~=$ \text{ }\text{ }\text{ }\text{ }\text{ }\text{ }\text{ }\text{ }\text{ }\text{ } \text{ } \text{ } \text{ } \text{ } \text{ } \text{ } \text{ } \text{ } \text{ } \text{ } \text{ } \text{ } \text{ } $(n-4) \cdot \binom{n-|V(G)|}{2};$\\
&  & &  \vspace{0.1cm} $\beta_{n-2, n + \alpha+2\beta-3}(R/I) = n-2.$\\ \hline

  &\vspace{0.1cm}$\girth(G) = 3$ &  \vspace{0.1cm}$G$ does not &
\vspace{0.1cm}$\beta_{n-2, n + 3\alpha-3}(R/I) =c_3(G)$.\\
\cline{2-2}\cline{4-4} $(n-2,a_2+2)$ & $\girth(G) \neq 3$  &satisfy both $(G_1)$ and $(G_2)$ & $\beta_{n-2,n+2\alpha+\beta-3}(R/I)= \underset{1\leq i \leq n}\sum \binom{\deg(i)}{2}.$\\

\hline
\end{tabular}
\end{center}
\end{Theorem}

\newpage
\begin{Theorem} \label{extremal5}
Let $(\alpha, \beta)=(4, 2)$. The extremal Betti numbers of $R/I$ are:
\begin{center}
\begin{tabular}{ |p{2.5cm}|p{2.5cm}|p{3cm}|p{6.5cm}| }
\hline
Corners of Betti diagram & \multicolumn{2}{c|}{Structure of $G$} & Extremal Betti numbers of $R/I$ \\ \hline 
 &  & $G$ does not satisfy $(G_3)$ & 
$\beta_{n-1,n+6}(R/I)= 4p_0(G)+p_1(G);$\\
& && $\beta_{n-2, n + 9}(R/I) =c_3(G)$.\\

\cline{3-4} &$\girth(G) = 3$  &$G$ satisfies $(G_3),$ not $(G_4)$ &$\beta_{n-1,n+5}(R/I)=r(G);$ \\
&  & &   $\beta_{n-2, n + 9} (R/I)=c_3(G)$.\\

\cline{3-4} & & $G$ satisfies both of $(G_3)$, $(G_4)$ &$\beta_{n-1,n+3}(R/I)= \binom{n-|V(G)|}{3};$\\
&  & but not $(G_5)$ &   $\beta_{n-2, n + 9}(R/I) =c_3(G)$.\\

\cline{2-4}$(n-1,a_1+1);$ & &   $G$ satisfies $(G_1)$& $\beta_{n-1,n+6}(R/I)= 4p_0(G)+p_1(G);$\\
$(n-2,a_2+2)$ & &(so $G$ does not satisfy $(G_3)$) & $\beta_{n-2, n + 7}(R/I) =\underset{1\leq i \leq n}\sum \binom{\deg(i)}{2}.$\\

\cline{3-4} &$\girth(G)\neq 3,$ $G$ contains a vertex of degree $\geq 2$ & $(G)$ does not satisfy $(G_1)$ and $|V(G)|\leq n-3$ &
\vspace*{0.1cm}$\beta_{n-1,n+3}(R/I)= \binom{n-|V(G)|}{3};$ \\
& & (so $G$ satisfies both $(G_3),(G_4),$ not $(G_5)$) &  $\beta_{n-2, n +7} (R/I)=\underset{1\leq i \leq n}\sum \binom{\deg(i)}{2}.$\\

\cline{2-4} &$G$ consists of  & $|E(G)|= 1$ &$\beta_{n-1,n+3}(R/I)= \binom{n-2}{3};$\\
& disjoint edges & &  $\beta_{n-2, n + 5}(R/I) = n-2.$\\ \hline

\vspace*{0.1cm} $(n-1,a_1+1)$& $G$ consists of disjoint edges & \vspace*{0.1cm}$|E(G)|\neq 1$ & \vspace*{0.1cm}$\beta_{n-1,n+6}(R/I)~=~ 4p_0(G)$\\
&&&$\hspace*{2.8cm} = 4\binom{|E(G)|}{2}.$ \\

\hline

\vspace*{0.7cm}$(n-2,a_2+2)$ & \vspace*{0.1cm}$\girth(G) = 3$ & $G$ satisfies $(G_3),$ $(G_4),$ $(G_5)$ &\vspace*{0.1cm} $\beta_{n-2, n + 9}(R/I) =c_3(G)$.\\
\cline{2-4} & \vspace*{0.3cm}$\girth(G) \neq 3$& $(G)$ does not satisfy $(G_1),$ and $|V(G)|\geq n-2$ & \vspace*{0.1cm}$\beta_{n-2,n+7}(R/I)= \underset{1\leq i \leq n}\sum \binom{\deg(i)}{2}.$\\
\hline
\end{tabular}
\end{center}
where $r(G)  = \dfrac{1}{3}\underset{\underset{i\in V(G)\backslash (N(p)\cup N(q))}{pq\in G,} }\sum |\bigl((N(p)\cap N(q)\bigr)\setminus N(i)| $\\
\hspace*{2.2cm}$+ \dfrac{1}{2}\underset{\underset{i\in V(G)\backslash (N(p)\cup N(q))}{pq\in G,} }\sum |\bigl((N(p) \backslash N[q]) \cup (N(q)\backslash N[p])\bigr) \setminus N(i)| $\\

  \hspace*{2cm}$+ \underset{\underset{i\in V(G)\backslash (N(p)\cup N(q))}{pq\in G,} }\sum |[n] \setminus \bigl(N(p)\cup N(q)\cup N[i]\bigr) |.$
\end{Theorem}

We now describe the idea of the proofs of Theorems~\ref{extremal1}--\ref{extremal5}. For convenience, we write $a_i = a_i(R/I)$ for $i = 1, 2$.  As it was mentioned in Introduction, $R/I$ has at most two extremal Betti numbers, which are $\beta_{n-1, n + a_1}(R/I)$ and $\beta_{n-2, n + a_2}(R/I)$. Moreover, $R/I$ has two extremal Betti numbers if and only if it is not Cohen-Macaulay and $a_1 + 1 < a_2 + 2$. The Cohen-Macaulayness of $R/I$ was studied in \cite{MN}, while its $a_i$-invariants were determined in \cite{MT}. They were expressed in terms of $\alpha, \beta$ and the structure of $G$.

On the other hand, $\beta_{n-2, n + a_2}(R/I)$ is the unique extremal Betti number of $R/I$ if and only if $R/I$ is Cohen-Macaulay, while $\beta_{n-1, n + a_1}(R/I)$ is the unique extremal Betti number only when $a_1 + 1 \geq a_2 + 2$.

It follows that we need to compute $\beta_{n-1, n + a_1}(R/I)$ and $\beta_{n-2, n + a_2}(R/I)$ when they are extremal. In view of Theorem~\ref{corner}, to compute the extremal $\beta_{n-1, n + a_1}(R/I)$ amounts to calculate the number of vectors $\b \in \mathbb Z^n$ such that $|\b| = a_1$ and $H^1_{\mm}(R/I)_{\b} \neq 0$, and for each such $\b$, compute the dimension 
$$ 
\dim_K H^1_{\mm}(R/I)_{\b} = \dim_K\widetilde{H}_{-|\CS_\b|}(\Delta_\b(I); K)
$$ 
in view of Takayama's formula. The same process applies for $\beta_{n-2, n + a_2}(R/I)$; see Proposition~\ref{3} below.

We give a complete proof of Theorem~\ref{extremal2} in Proposition~\ref{3} and Proposition~\ref{thm2}. Although tempting, we refrain from presenting the proofs of the other main theorems, which are similar to that of Theorem~\ref{extremal2} and would make the article's length increases significantly.


\section{The values of $\beta_{n-2,n+a_2}(R/I)$ when it is extremal}

The computation of the extremal $\beta_{n-2,n+a_2}(R/I)$ is somewhat simpler than that of $\beta_{n-1, n+a_1}(R/I)$. We have the following result.
 
\begin{Proposition}{\label{3}} 
If $\beta_{n-2,n+a_2}(R/I)$ is an extremal Betti number of $R/I,$ then 
\begin{equation*}
\beta_{n-2,n+a_2}(R/I)=
\begin{cases}
(n-2)|E(G)| & \text{ if $G$ consists of disjoint edges;}\\
c_3(G) & \text{ if $\girth(G)=3$;}\\
a(G) + (\alpha - 2)c_4(G) & \text{ if $\beta =1$; $\girth(G)\neq 3$}\\
&\text{ and $G$ contains a vertex of degree $\geq 2$;}\\
\underset{1\leq i \leq n}\sum \binom{\deg(i)}{2}   & \text{ if $\beta \neq 1$; $\girth(G)\neq 3$}\\
&\text{ and $G$ contains a vertex of degree $\geq 2$.}\\
\end{cases} 
\end{equation*}
  \end{Proposition}

 \begin{proof}
We know that $\beta_{n-2, n + a_2}(R/I)$ is an extremal Betti number if and only if  $(2, a_2 + 2)$ is a corner of the local cohomology diagram of $R/I$, and 
$$
\beta_{n-2, n + a_2}(R/I) = \dim_K H^2_{\mm}(R/I)_{a_2},
$$ 
by Theorem~\ref{corner}. Let $\b = (b_1,b_2,\ldots,b_n) \in\mathbb{Z}^n$ be such that $|\b| = a_2$ and $H^2_{\mm}(R/I)_{\b} \neq 0$. By Takayama's formula, we have
$$ 
0 \neq \dim_K H^2_{\mm}(R/I)_{\b} = \dim_K\widetilde{H}_{1-|\CS_{\b}|}(\Delta_\b(I);K).
$$ 
We need to calculate the number of such vectors $\b$, and for each $\b$, we determine $\dim_K\widetilde{H}_{1 - |\CS_{\b}|}(\Delta_\b(I);K)$.

By \cite[Lemmas 3.1--3.4]{MT}, $H^2_{\mm}(R/I)_{\b}\neq 0$ if and only if $\b\in \mathbb{N}^n$ and $\Delta_\b(I)$ must contain a cycle, say $\{1,2\},\{2,3\},\ldots,\{t-1,t\},\{1,t\}$ for some $t\geq 3$. By Lemma \ref{L2},  
$\sigma_{i,i+1}(\b)\leq w_{i,i+1}-1$ for every $i=1,\ldots,t-1$ and $\sigma_{1,t}(\b)\leq w_{1,t}-1$. Hence, \begin{equation}\label{Ineq}
 (t-2)|\b| \leq t|\b|-2(b_1+b_2+\cdots+b_t)\leq \sum^{t-1}_{i=1}w_{i,i+1}+w_{1,t}-t \leq t\alpha - t.
\end{equation}
It follows that \begin{equation}\label{Ineq2}
|\b| \le \dfrac{t\alpha-t}{t-2}=\alpha -1+\dfrac{2\alpha-2}{t-2} \leq 3\alpha -3.
\end{equation}
We distinguish three cases as follows.

\medskip

\noindent{\bf Case 1. $\girth(G)=3$.} Then using \cite[Theorem 3.5]{MT}, $|\b|=a_2(R/I)=3\alpha-3.$ It follows from (\ref{Ineq}) and (\ref{Ineq2}) that 
\begin{equation*}
\begin{cases}
t=3;\\
w_{1,2}=w_{1,3}=w_{2,3}=\alpha;\\
\sigma_{1,2}(\b)=\sigma_{1,3}(\b)=\sigma_{2,3}(\b)=\alpha-1;\\
b_i=0 \text{ for all } i\geq 4.
\end{cases}
\end{equation*}
Therefore, $b_1=b_2=b_3=\alpha-1$ and $b_i=0 \text{ for all } i\geq 4.$ By Lemma \ref{L2}, we have 
$$\Delta_\b(I)=\{\{1,2\},\{1,3\},\{2,3\}\}.$$
Hence $\dim_{K}H^2_{\mm}(R/I)_{\b}=\dim_K\widetilde{H}_{1}(\Delta_\b(I);K)=1.$ 

On the other hand, each 3-cycle of $G$ corresponds to such a vector $\b$. Thus, in Case 1 we get $$\beta_{n-2,n+a_2}(R/I)= c_3(G).$$

\medskip

\noindent{\bf Case 2. $\girth(G)\neq 3$ and $G$ contains a vertex of degree $\geq 2.$} By \cite[Theorem 3.5]{MT}, $|\b|=a_2=2\alpha+\beta - 3.$ We consider two subcases.

\medskip

\noindent {\bf Subcase 2.1.} $t=3$. Then $|\b|\leq w_{1,2}+w_{2,3}+w_{1,3} - 3 \leq 2\alpha +\beta -3$ by (\ref{Ineq}) and the fact that $\girth G \neq 3$. Therefore, there exist $i\neq j \in \{1,2,3\}$ such that $w_{i,j}=\beta$, say $i=1, j=2.$ Then $|\b|=2\alpha+\beta - 3$ if and only if 
\begin{equation*}
\begin{cases}
w_{1,2}=\beta, w_{1,3}=w_{2,3}=\alpha;\\
\sigma_{1,2}(\b)=\beta -1,\sigma_{1,3}(\b)=\sigma_{2,3}(\b)=\alpha-1;\\
b_i=0 \text{ for all } i\geq 4.
\end{cases}
\end{equation*}
Hence, $b_1=b_2=\alpha -1, b_3=\beta -1$ and $b_i=0 \text{ for all } i\geq 4,$ i.e., 
$$
\b=(\alpha-1,\alpha-1,\beta-1,0,\dots,0).
$$ 
Using Lemma \ref{L2}, $\Delta_\b(I)$ is $\{1,2\}\cup\{\{1,i\}\mid 1i\in G\}\cup \{\{2,i\}\mid 2i\in G\} $ if $\beta=1$; or $\{\{1,2\},\{1,3\},\{2,3\}\}$ otherwise. By Lemma \ref{Euler}, we have
\begin{equation*}
\dim_{K}H^2_{\mm}(R/I)_{\b}=\dim_K\widetilde{H}_{1}(\Delta_\b(I);K)=
\begin{cases}
|N(1)\cap N(2)| & \text{ if $\beta =1$;}\\
1 &\text{ otherwise.}
\end{cases}
\end{equation*}
Moreover, each $ij\not\in G$ which satisfies $N(i)\cap N(j)\neq \emptyset$ corresponds to such a vector $\b$ if $\beta=1$, and each pair of joint edges of $G$ corresponds to such a  vector $\b$ if $\beta\neq 1.$

\medskip

\noindent {\bf Subcase 2.2.} $t\neq 3$. Using (\ref{Ineq}), we have 
$$
|\b|\le \dfrac{t\alpha-t}{t-2}=\alpha -1+\dfrac{2\alpha-2}{t-2}\le \alpha -1+\dfrac{2\alpha-2}{2}\leq 2\alpha+\beta-3.
$$
So $|\b|=2\alpha+\beta -3$ if and only if 
\begin{equation*}
\begin{cases}
t=4;\\
\beta=1;\\
w_{1,2}=w_{2,3}=w_{3,4}=w_{1,4}=\alpha;\\
\sigma_{1,2}(\b)=\sigma_{2,3}(\b)=\sigma_{3,4}(\b)=\sigma_{1,4}(\b)=\alpha-1;\\
b_i=0 \text{ for all } i\geq 5.
\end{cases}
\end{equation*}
Therefore, $b_1=b_3, b_2=b_4, b_1+b_2=\alpha-1$ and $b_i=0 \text{ for all } i\geq 5,$ i.e., 
$$
\b=(b_1,\alpha-1-b_1,b_1,\alpha-1-b_1,0,\dots,0) \text{ for } 0\leq b_1\leq \alpha-1.
$$

If $0<b_1<\alpha-1,$ then 
$$
\Delta_\b(I) = \{\{1,2\},\{2,3\},\{3,4\},\{1,4\}\}
$$ 
as $\beta =1$, so that $\dim_{K}H^2_{\mm}(R/I)_{\b}=\dim_K\widetilde{H}_{1}(\Delta_\b(I);K) = 1$. Conversely, each $4$-cycle of $G$ and each $0 < b_1  < \alpha - 1$ correspond to such a vector $\b$.

If $b_1=0$ or $b_1=\alpha-1$, then $\D_\b(I)$ contains a 3-cycles and $\b$ was constructed in {Subcase 2.1}.

Therefore, in Case 2 we obtain 
\begin{equation*}
\beta_{n-2,n+a_2}(R/I)=
\begin{cases}
\underset{ij\not\in G}\sum|N(i)\cap N(j)| + (\alpha-2)c_4(G) & \text{ if $\beta =1$;}\\
\underset{1\leq i \leq n}\sum \binom{\deg(i)}{2} &\text{ otherwise.}
\end{cases}
\end{equation*}

\noindent{\bf Case 3.} {\bf Every vertex of $G$ has degree $1$.} Using \cite[Theorem 3.5]{MT}, we have $|\b|=a_2(R/I)=\alpha+2\beta - 3.$ On the other hand, if $t>3$ and $t$ is even, then  $(t-2)|\b|\leq\sum^{t-1}_{i=1}w_{i,i+1}+w_{1,t}-t\leq \dfrac{t}{2}\alpha+\dfrac{t}{2}\beta-t$ by (\ref{Ineq}). Hence,
$$|\b| \leq \dfrac{\alpha+\beta-2}{2}+\dfrac{\alpha+\beta-2}{t-2} \le \alpha + \beta -2 \le\alpha+2\beta -3.$$
So $|\b|=\alpha+2\beta -3$ if and only if 
\begin{equation*}
\begin{cases}
t=4;\\
\beta=1;\\
w_{1,2}=w_{3,4}=\alpha, w_{2,3}=w_{1,4}=\beta;\\
\sigma_{1,2}(\b)=\sigma_{3,4}(\b)=\alpha-1,\sigma_{2,3}(\b)=\sigma_{1,4}(\b)=\beta-1;\\
b_i=0 \text{ for all } i\geq 5,
\end{cases}
\end{equation*}
or 
\begin{equation*}
\begin{cases}
t=4;\\
\beta=1;\\
w_{1,2}=w_{3,4}=\beta, w_{2,3}=w_{1,4}=\alpha;\\
\sigma_{1,2}(\b)=\sigma_{3,4}(\b)=\beta-1,\sigma_{2,3}(\b)=\sigma_{1,4}(\b)=\alpha-1;\\
b_i=0 \text{ for all } i\geq 5.
\end{cases}
\end{equation*}
Hence, $2(\alpha-1)=\sigma_{1,2}(\b)+\sigma_{3,4}(\b)=\sigma_{2,3}(\b)+\sigma_{1,4}(\b)=2(\beta-1),$ a contradiction.

If $t>3$ and $t$ is odd, then $(t-2)|\b|\leq\sum^{t-1}_{i=1}w_{i,i+1}+w_{1,t}-t\le \dfrac{t-1}{2}\alpha+\dfrac{t+1}{2}\beta-t$ by (\ref{Ineq}). Therefore,
$$
|\b| \le \dfrac{\alpha+\beta-2}{2}+\dfrac{\alpha+3\beta-4}{2(t-2)}\le \dfrac{\alpha+\beta-2}{2} + \dfrac{\alpha+3\beta-4}{6}< \alpha+2\beta -3,
$$
another contradiction.

If $t=3$, then since every vertex of $G$ has degree $1$ and by (\ref{Ineq}), $$|\b|\leq w_{1,2}+w_{2,3}+w_{1,3}-3\leq \alpha+2\beta-3.$$
So $|\b|=\alpha+2\beta -3$ if and only if  
\begin{equation*}
\begin{cases}
w_{i,j}=\alpha, w_{i,k}=w_{j,k}=\beta;\\
\sigma_{i,j}(\b)=\alpha-1, \sigma_{i,k}(\b)=\sigma_{j,k}(\b)=\beta-1;\\
b_l=0 \text{ for all } l\geq 4,
\end{cases}
\end{equation*}
for $\{i,j,k\}=\{1,2,3\}$, say $i=1,j=2$. Then $\b=(\beta -1,\beta-1,\alpha-1,0,\dots,0).$ 
It follows that $\D_\b(I)$ is $\{\{1,2\},\{1,3\},\{2,3\}\}$ if $2\beta-1>\alpha$, or $\{\{1,2\},\{1,3\},\{2,3\}\}\cup \{\{3,i\}\mid 3i\in G\}$ if $\beta \neq 1$ and $2\beta -1 \leq \alpha$, or $\{\{i,j\}\mid ij\in G\}\cup \{\{3,i\}\mid 1\leq i\neq 3\leq n\}$ otherwise. Therefore,
\begin{equation*}
\dim_{K}H^2_{\mm}(R/I)_{\b}=\dim_K\widetilde{H}_{1}(\Delta_\b(I);K)=
\begin{cases}
1 & \text{ if $\beta \neq 1$;} \\
|E(G)|&\text{ if $\beta =1$ and $3\not\in V(G)$;} \\
|E(G)|-1&\text{ if $\beta =1$ and $3\in V(G)$.}
\end{cases} 
\end{equation*}
On the other hand, the number of vectors $\b$ is $(n-2)|E(G)|$ if $\beta\neq 1$; or $n-2$ if $\beta=1$ and $|E(G)|=1$; or $n$ if $\beta=1$ and $|E(G)|\geq 2.$ Thus, in this case we obtain
\begin{align*}
\beta_{n-2,n+a_2}(R/I) = 
\begin{cases}
(n-2) |E(G)| & \text{ if $\beta \neq 1$;}\\
(n-2) |E(G)| &\text{ if $\beta =1$, $|E(G)|=1$;}\\
(n-2|E(G)|).|E(G)| + 2|E(G)|.(|E(G)|-1) &\text{ if $\beta =1$, $|E(G)|\geq 2$;}
\end{cases}
\end{align*} 
i.e., $\beta_{n-2,n+a_2}(R/I)=(n-2) |E(G)|.$

From all the above cases, the proof is completed.
\end{proof}

\section{The values of $\beta_{n-1, n+a_1}(R/I)$ when it is extremal and $\alpha\geq \beta +3$}

As we may observe from the previous section, to count the numbers of vectors $\b$ such that $H^2_{\mm}(R/I)_{\b}\neq 0$ or $H^1_{\mm}(R/I)_{\b}\neq 0$ is amount to solve a system of linear inequalities in the integers. Therefore, we start this section with two auxiliary computations.

\begin{Lemma} \label{sol1} 
Consider the system of linear inequalities
\begin{equation} \label{system1}
\begin{cases}
0 < x < \alpha - 1, \; 0 < y < \alpha - 1 \\
\beta \leq x + y \leq 2\alpha - \beta - 2 \\
\beta + 1 - \alpha \leq x - y \leq \alpha - \beta - 1
\end{cases}.
\end{equation}
The number of integer solutions $(x, y)$ of the system \eqref{system1} is
 \begin{equation*}
     \begin{cases}
     (\alpha - \beta)^2 + (\alpha - \beta - 1)^2 - \overline{\alpha - 2\beta + 1}^2 & \text{ if $\alpha$ is odd;} \\
         2(\alpha - \beta) (\alpha - \beta - 1)- \overline{\alpha - 2\beta}\,(\alpha - 2\beta + 2) & \text{ if $\alpha$ is even.}
     \end{cases}
 \end{equation*}
\end{Lemma}

\begin{proof}
We assume that $\alpha$ is odd; similarly for the case $\alpha$ is even. First, let $S$ be the set of integer solutions of the system
\begin{equation} \label{system2} 
\begin{cases}
\beta \leq x + y \leq 2\alpha - \beta - 2 \\
\beta + 1 - \alpha \leq x - y \leq \alpha - \beta - 1
\end{cases}.
\end{equation}
We calculate $|S|$. Let $x + y = \beta + k$ with $k \in \mathbb Z$. Then $y = (\beta + k) - x$, and  
\begin{align*} 
\eqref{system2} \Leftrightarrow \begin{cases}
\beta \leq \beta + k \leq 2\alpha - \beta - 2 \\
\beta + 1 - \alpha \leq x - y \leq \alpha - \beta - 1
\end{cases} \Leftrightarrow 
\begin{cases}
0 \leq k \leq 2\alpha - 2\beta - 2 \\
\dfrac{2\beta - \alpha + k + 1}{2} \leq x \leq \dfrac{\alpha+k-1}{2} 
\end{cases}.
\end{align*}
For a fixed $k$, the number of integers $x$ satisfying the above linear inequalities equals
$$
\begin{cases}
\alpha - \beta & \text{ if $k$ is even;} \\
\alpha - \beta - 1 & \text{ if $k$ is odd.}
\end{cases}
$$
Under the condition $0 \leq k \leq 2\alpha - 2\beta - 2$, there are $(\alpha - \beta)$ even integers $k$ and $(\alpha - \beta - 1)$ odd integers $k$. It follows that 
$$
|S| = (\alpha - \beta)^2 + (\alpha - \beta - 1)^2.
$$
Next, we count the number of integer solutions $(x, y)$ of \eqref{system2} which are not solutions of \eqref{system1}. Such a pair $(x, y)$ belongs to one of the following set: 
\begin{align*}
S_1 = \{(x, y) \in S \mid x \geq \alpha - 1\}, \quad S_2 = \{(x, y) \in S \mid x \leq 0\}, \\ 
S_3 = \{(x, y) \in S \mid y \geq \alpha - 1\}, \quad S_4 = \{(x, y) \in S \mid y \leq 0\}.
\end{align*}
Observe that the sets $S_i$'s are pairwise disjoint. Indeed, if $(x, y) \in S_1$, then 
$$
y \leq (2\alpha - \beta - 2) - x \leq (2\alpha - \beta - 2) - (\alpha - 1) = \alpha - \beta - 1 < \alpha - 1,
$$
and
$$
y \geq x - (\alpha - \beta - 1) \geq (\alpha - 1) - (\alpha - \beta - 1) = \beta > 0.
$$
This shows that $S_1$ is disjoint from $S_2 \cup S_3 \cup S_4$. Similar arguments confirm the above observation.

Let $(x, y) \in S_1$ and again write $x + y = \beta + k$ for $k \in \mathbb{Z}$ and $0 \leq k \leq 2\alpha - 2\beta - 2$. Since $(x, y)$ satisfies \eqref{system2}, we have $\dfrac{2\beta - \alpha + k + 1}{2} \leq x \leq \dfrac{\alpha+k-1}{2}$. However, $\dfrac{2\beta - \alpha + k + 1}{2} \leq \dfrac{\alpha - 1}{2}  < \alpha - 1$ as $k \leq 2\alpha - 2\beta - 2$. On the other hand, $\alpha - 1 \leq \dfrac{\alpha +k - 1}{2}$ if and only if $k \geq \alpha - 1$. Thus, for each $k = \alpha - 1, \alpha, \ldots, 2\alpha - 2\beta - 2$ (assuming $2\beta + 1 \leq \alpha$ so that $\alpha - 1 \leq 2\alpha - 2\beta - 2$), we count the number of integers $x$ satisfying
$$
\alpha - 1 \leq x \leq \dfrac{\alpha + k - 1}{2}.
$$
It follows that (noting that $\alpha$ is odd)
$$
|S_1| = 2\biggr(1 + 2 + \ldots + \dfrac{\alpha - 2\beta -1}{2}\biggl) + \dfrac{\alpha - 2\beta +1}{2}  = \biggr(\dfrac{\alpha - 2\beta + 1}{2}\biggl)^2.
$$
Analogously, we find that $|S_2| = |S_3| = |S_4| = |S_1|$. Consequently, the number of integer solutions of \eqref{system1} equals
\begin{align*}
& |S| - (|S_1| + |S_2| + |S_3| + |S_4|) \\
& = \begin{cases}
(\alpha - \beta)^2 + (\alpha - \beta - 1)^2 - (\alpha - 2\beta + 1)^2 & \text{ if } \alpha \geq 2\beta + 1; \\
(\alpha - \beta)^2 + (\alpha - \beta - 1)^2 & \text{ if } \alpha < 2\beta + 1.
\end{cases}
\end{align*}
\end{proof}

\begin{Lemma} \label{sol2}
Let $\beta > 1$ and $\alpha \geq \beta +2$. Consider the system of linear inequalities
\begin{equation} \label{system3}
\begin{cases}
0 < x < \alpha - 1, \; 0 < y < \beta - 1 \\
\beta \leq x + y \leq \alpha - 2 \\
1 \leq x - y \leq \alpha - \beta - 1
\end{cases}.
\end{equation}
The number of integer solutions $(x, y)$ of the system \eqref{system3} equals
\begin{equation*}
\begin{cases}
\dfrac{1}{2}\biggl[(\alpha - \beta - 1)^2 - \overline{\alpha - 2\beta + 1}^2\biggr] & \text{ if $\alpha$ is odd, $\beta$ is even;} \\
\dfrac{1}{4}\biggl[(\alpha - \beta)^2 + (\alpha - \beta - 2)^2 - 2\, \overline{\alpha - 2\beta + 1}^2\biggr] & \text{ if $\alpha, \beta$ are odd;} \\
\dfrac{1}{2}\biggl[{(\alpha - \beta)(\alpha - \beta - 2)}- \overline{\alpha - 2\beta}\, (\alpha - 2\beta + 2)\biggr] & \text{ if $\alpha, \beta$ are even;} \\
\dfrac{1}{2}\biggl[(\alpha - \beta - 1)^2 - \overline{\alpha - 2\beta}\, (\alpha - 2\beta + 2)\biggr] & \text{ if $\alpha$ is even, $\beta$ is odd.}
\end{cases}
\end{equation*}
\end{Lemma}

\begin{proof}
We consider the case when $\alpha$ is odd and $\beta$ is even, while the other cases are similar. Let $T$ be the set of integer solutions of the system
\begin{equation} \label{system4} 
\begin{cases}
\beta \leq x + y \leq \alpha - 2 \\
1 \leq x - y \leq \alpha - \beta - 1 
\end{cases}.
\end{equation}
We calculate $|T|$. Let $x + y = \beta + k$ with $k \in \mathbb Z$. Then $y = (\beta + k) - x$, and  
\begin{align*} 
\eqref{system4} \Leftrightarrow \begin{cases}
\beta \leq \beta + k \leq \alpha - 2 \\
1 \leq x - y \leq \alpha - \beta - 1
\end{cases} \Leftrightarrow 
\begin{cases}
0 \leq k \leq \alpha - \beta - 2 \\
\dfrac{\beta+k+1}{2} \leq x \leq \dfrac{\alpha+k-1}{2} 
\end{cases}.
\end{align*}
For a fixed $k$ (even or odd), the number of integers $x$ satisfying the above linear inequalities equals
$\dfrac{\alpha - \beta - 1}{2}$. The number of integers $k$ satisfying $0 \leq k \leq \alpha - \beta - 2$  is $\alpha - \beta - 1$. It follows that 
$$
|T| = \dfrac{(\alpha - \beta - 1)^2}{2}.
$$
Next, we count the number of integer solutions $(x, y)$ of \eqref{system4} which are not solutions of \eqref{system3}. Such a pair $(x, y)$ belongs to one of the following set: 
\begin{align*}
T_1 = \{(x, y) \in T \mid x \geq \alpha - 1\}, \quad T_2 = \{(x, y) \in T \mid x \leq 0\}, \\ 
T_3 = \{(x, y) \in T \mid y \geq \beta -1\}, \quad T_4 = \{(x, y) \in T \mid y \leq 0\}.
\end{align*}
Observe that if $(x,y)\in T,$ then $\beta + 1 \leq (x-y) + (x+ y)\leq (\alpha -2) + (\alpha -\beta -1),$ hence $\beta +1 \leq 2x \leq 2\alpha -\beta -3 \leq 2\alpha -4.$ It shows that $T_ 1 = T_2 = \emptyset$. Clearly $T_3$ and $T_4$ are disjoint, and  for $(x, y) \in T_3 \cup T_4$, we deduce that 
$$ 
\beta \leq x \leq \alpha - \beta - 1.
$$  
Hence $T_3 \cup T_4 = \emptyset$ when $\alpha < 2\beta + 1$. 
 
Assume that $\alpha \geq 2\beta + 1$. Let $(x, y) \in T_3$.  We write again $x + y = \beta + k$ for $k \in \mathbb{Z}$. Then $x = (\beta + k) - y$, and  
\begin{align*} 
\eqref{system4} \Leftrightarrow 
\begin{cases}
0 \leq k \leq \alpha - \beta - 2 \\
\dfrac{2\beta - \alpha + k + 1}{2} \leq y \leq \dfrac{\beta+k-1}{2} 
\end{cases}.
\end{align*}
Clearly $\dfrac{2\beta - \alpha + k + 1}{2} < \beta - 1$ as $k \leq \alpha - \beta - 2$ and $\beta > 1$. On the other hand, $\beta - 1 \leq \dfrac{\beta + k - 1}{2}$ if and only if $k \geq \beta - 1$. Thus, for each $k = \beta - 1, \beta, \ldots, \alpha - \beta - 2$, we count the number of integers $y$ satisfying
$$
\beta - 1 \leq y \leq \dfrac{\beta + k - 1}{2}.
$$
It follows that (noting that $\alpha$ is odd)
\begin{equation*}
|T_3| =
 2\biggr(1 + 2 + \ldots + \dfrac{\alpha - 2\beta - 1}{2}\biggl) + \dfrac{\alpha - 2\beta + 1}{2}  = \dfrac{(\alpha - 2\beta + 1)^2}{4} .
\end{equation*}
In a similar manner, we obtain
$$
|T_4| = |T_3| = \dfrac{(\alpha - 2\beta + 1)^2}{4}.
$$
Consequently, the number of integer solutions $(x, y)$ of the system \eqref{system3} is
\begin{equation*}
|T| - |T_3 \cup T_4| = 
\begin{cases}
\dfrac{(\alpha - \beta - 1)^2}{2} - \dfrac{(\alpha - 2\beta + 1)^2}{2} & \text{ if } \alpha \geq 2\beta + 1; \\
\dfrac{(\alpha - \beta - 1)^2}{2} & \text{ if } \alpha < 2\beta + 1
\end{cases}
\end{equation*}
when $\alpha$ is odd and $\beta$ is even.
\end{proof}

We are now in a position to compute the values of the extremal $\beta_{n-1, n + a_1}(R/I)$ when $\alpha \geq \beta + 3$.

\begin{Proposition} \label{thm2}
Let $\alpha \geq \beta + 3$. If $\beta_{n-1, n + a_1}(R/I)$ is an extremal Betti number of $R/I$, then 
\begin{equation*}
\beta_{n-1, n + a_1}(R/I) = \begin{cases} 
B_1(\alpha, \beta) & \text{ if $G$ satisfies $(G_1)$}; \\
B_2(\alpha, \beta) & \text{ if $G$ does not satisfy $(G_1)$}.
\end{cases}
\end{equation*}
\end{Proposition}

\begin{proof}
We know that $\beta_{n-1, n + a_1}(R/I)$ is an extremal Betti number of $R/I$ if and only if $H^1_{\mm}(R/I)_{a_1}$ is an extremal local cohomology, and
$$
0 \neq \beta_{n-1, n + a_1}(R/I) = \dim_K H^1_{\mm}(R/I)_{a_1},
$$
by Theorem~\ref{corner}. In particular, $R/I$ is not Cohen-Macaulay. Let $\b = (b_1,b_2,\ldots,b_n) \in\mathbb{Z}^n$ be such that $|\b| = a_1(R/I)$ and $H^1_{\mm}(R/I)_{\b} \neq 0$. It follows from Lemma~\ref{L3} that $\b \in \mathbb{N}^n$ and $\D_\b(I)$ is a disconnected graph, and from Takayama's formula, 
$$ 
\dim_K H^1_{\mm}(R/I)_{\b} = \dim_K\widetilde{H}_{0}(\Delta_\b(I);K).
$$ 
We need to count the number of such vectors $\b$ and calculate $\dim_K\widetilde{H}_{0}(\Delta_\b(I);K)$ for each $\b$. To do so we base on the configurations of nonzero coordinates of $\b$. Assume that $\{1, 2\}, \{3, 4\}$ belong to different connected components of $\D_\b(I)$. We distinguish two cases as follows.

\medskip

\noindent {\bf Case 1. $G$ satisfies $(G_1)$}. By \cite[Theorem 2.2]{MT}, $|\b| = 2\alpha - 2$. Using Lemma~\ref{L2}, we have $\sigma_{1,2}(\b) \leq w_{1,2} - 1 \leq \alpha - 1$ and $\sigma_{3,4}(\b) \leq w_{3,4} - 1 \leq \alpha -1$. Hence
$$
2\alpha - 2 = |\b| \leq |\b| + \sum_{i \geq 5} b_i  = \sigma_{1,2}(\b) + \sigma_{3,4}(\b) \leq 2\alpha - 2.
$$
Consequently, we get $\sigma_{1,2}(\b) = \sigma_{3,4}(\b) = \alpha - 1$ and $b_i = 0$ for all  $i \geq 5$. It follows that $b_1 + b_2 = b_3 + b_4 = \alpha - 1$, so the vector $\b$ has the form
\begin{equation} \label{form1}
(b_1, \alpha - 1 - b_1, b_3, \alpha - 1 - b_3, 0, \ldots, 0). 
\end{equation}
We note that $w_{1,2} = w_{3,4} = \alpha$, that is, $12, 34 \in G$. Assume that $\alpha$ is odd and $\beta$ is even (the other cases are analogous).

\medskip

\noindent {\bf Subcase 1.1.} The induced subgraph $G[1,2,3,4]$ is $\{12, 34\}$, i.e., $(12, 34)$ is a disconnected pair of edges of $G$. By Lemma~\ref{L2}, we have
$$
\sigma_{1,3}(\b)\geq w_{1,3} = \beta, \sigma_{1,4}(\b) \geq w_{1,4} = \beta, \sigma_{2,3}(\b) \geq w_{2,3} = \beta, \text{ and } \sigma_{2,4}(\b) \geq w_{2,4} = \beta.
$$
It follows that
\begin{equation} \label{eq1}
\begin{cases}
\beta \leq b_1 + b_3 \leq 2\alpha - \beta - 2 \\
\beta + 1 - \alpha \leq b_1 - b_3 \leq \alpha - \beta - 1.
\end{cases}
\end{equation}
Thus, the number of vectors $\b$ of the form~\eqref{form1} in Subcase 1.1 equals the number of integral points $(b_1, b_3)$ in the square given by the system of linear inequalities~\eqref{eq1} such that $\b \in \mathbb{N}^n$. By Lemma~\ref{sol1}, the number of vectors $\b$ with $0 < b_1, b_3 < \alpha - 1$ and $b_1, b_3$ satisfying \eqref{eq1} equals
$$
(\alpha - \beta)^2 + (\alpha - \beta - 1)^2 - \overline{\alpha - 2\beta + 1}^2.
$$
For such a vector $\b$, we get $\mathcal F(\Delta_\b(I))=\{\{1,2\}, \{3,4\}\}$ by Lemma~\ref{L2}, so 
$$
\dim_K H^{1}_{\mm}(R/I)_\b = \dim_K\widetilde{H}_{0}(\Delta_\b(I);K) = 1.
$$ 
Conversely, for each disconnected pair of edges $(ij, pq)$ of $G$, there are
$$
(\alpha - \beta)^2 + (\alpha - \beta - 1)^2 - \overline{\alpha - 2\beta + 1}^2.
$$
vectors $\b$ of the form
$
(0, \ldots, \underbrace{b_1}_{i}, \ldots, \underbrace{\alpha-1-b_1}_{j}, \ldots, \underbrace{b_3}_{p}, \ldots, \underbrace{\alpha - 1 - b_3}_{q}, 0, \ldots, 0)
$
with $0 < b_1, b_3 < \alpha - 1$ and $b_1, b_3$ satisfying the system \eqref{eq1}.

If $b_1 = \alpha-1$ and $0 < b_3 < \alpha - 1$, then $\b = {\b'} = (\alpha-1, 0, b_3, \alpha - 1 - b_3, 0, \ldots, 0)$. The system \eqref{eq1} turns into
\begin{equation} \label{eq2}
\begin{cases}
\beta - \alpha + 1 \leq b_3 \leq \alpha - \beta - 1 \\
\beta \leq b_3 \leq 2\alpha - \beta - 2.
\end{cases}
\end{equation}
It follows that $\beta \leq b_3 \leq \alpha - \beta - 1$, so the number of vectors ${\b'}$ is $\overline{\alpha - 2\beta}$ (note that $\alpha$ is odd). By Lemma~\ref{L2}, $\F(\D_{\b'}(I)) = \{\{1,i\} \mid 1i \in G\} \cup \{\{3,4\}\}$, thus
$$
\dim_K H^{1}_{\mm}(R/I)_{\b'}=\dim_K\widetilde{H}_{0}(\Delta_{\b'}(I);K)= 1.
$$
The situation $b_1 = 0$ and $0 < b_3 < \alpha - 1$ is analogous. In general, we see that for each $pq \in G$ and $i\not\in N(p)\cup N(q)$, there are $\overline{\alpha - 2\beta}$ vectors ${\b'}$ as above with $\alpha-1$ in the $i$th coordinate and other two nonzero coordinates in the $p$th and $q$th.

If $b_1, b_3 \in \{0, \alpha-1\}$, then one of the two inequalities in the system \eqref{eq1} is not satisfied.

Consequently, in the Subcase 1.1 the number of vectors $\b$ in question equals
\begin{equation} \label{case1.1}
\biggl[(\alpha - \beta)^2 + (\alpha - \beta - 1)^2 - \overline{\alpha - 2\beta + 1}^2\biggr]p_0(G) + \overline{\alpha - 2\beta}\, b(G).
\end{equation}
For each vector $\b$, we have $\dim_K H^{1}_{\mm}(R/I^{(m)})_{\b} = 1$.

\medskip

\noindent {\bf Subcase 1.2.} The induced subgraph $G[1,2,3,4]$ is a path of length 3, i.e., $(12, 34)$ is a type-1 pair of edges. Without loss of generality, we may assume $G[1, 2, 3, 4] = \{12, 34, 24\}$.
By Lemma~\ref{L2}, we have
$$
\sigma_{2,4}(\b) \geq w_{2,4} = \alpha, \sigma_{1,3}(\b)\geq w_{1,3} = \beta, \sigma_{1,4}(\b) \geq w_{1,4} = \beta, \sigma_{2,3}(\b) \geq w_{2,3} = \beta.
$$
It follows that
\begin{equation} \label{eq3}
\begin{cases}
\alpha \leq b_1 + b_3 \leq 2\alpha - \beta - 2 \\
\beta + 1 - \alpha \leq b_1 - b_3 \leq \alpha - \beta - 1.
\end{cases}
\end{equation}

 If $0 < b_1, b_3 < \alpha - 1$ and $b_1, b_3$ satisfy \eqref{eq3}, then similarly to Lemma~\ref{sol1}, it can be computed that the number of vectors $\b$ of the form \eqref{form1} equals
$$
\dfrac{1}{2}\bigg[(\alpha - \beta - 1)(2\alpha - 2\beta - 1) - \overline{\alpha - 2\beta + 1}^2\biggr].
$$
For such a vector $\b$, we get $\mathcal F(\Delta_\b(I))=\{\{1,2\}, \{3,4\}\}$ by Lemma~\ref{L2}, so 
$$
\dim_K H^{1}_{\mm}(R/I)_\b = \dim_K\widetilde{H}_{0}(\Delta_\b(I);K) = 1.
$$ 
Conversely, for each type-1 pair of edges $(ij, pq)$ of $G$, there are
$$
\dfrac{1}{2}\bigg[(\alpha - \beta - 1)(2\alpha - 2\beta - 1) - \overline{\alpha - 2\beta + 1}^2\biggr]
$$
vectors $\b$ of the form
$
(0, \ldots, \underbrace{b_1}_{i}, \ldots, \underbrace{\alpha-1-b_1}_{j}, \ldots, \underbrace{b_3}_{p}, \ldots, \underbrace{\alpha - 1 - b_3}_{q}, 0, \ldots, 0)
$
with $0 < b_1, b_3 < \alpha - 1$ and $b_1, b_3$ satisfying the system \eqref{eq3}. We note that these vectors arise differently to those in the Subcase 1.1 which correspond to disconnected pairs of edges.

If either $b_1$ or $b_3$ equals $0$ or $\alpha-1$, then $\b$ has the form $\b'$ as in the Subcase 1.1 and these vectors are already counted. If both $b_1$ and $b_3$ belong to $\{0, \alpha-1\}$, then the first inequality in \eqref{eq3} is not satisfied. Therefore, in the Subcase 1.2 the number of vectors $\b$ is 
\begin{equation} \label{case1.2}
\dfrac{1}{2}\bigg[(\alpha - \beta - 1)(2\alpha - 2\beta - 1) - \overline{\alpha - 2\beta + 1}^2\biggr] p_1(G).
\end{equation}
Also, for each of these vectors $\b$, we have $\dim_K H^{1}_{\mm}(R/I^{(m)})_{\b} = 1$.

\medskip

\noindent {\bf Subcase 1.3.} The induced subgraph $G[1,2,3,4]$ is a path of length 4 which is not a 4-cycle, i.e., $(12, 34)$ is a type-2 pair of edges. We may assume that $G[1,2,3,4] = \{12, 34, 14, 24\}$. By Lemma~\ref{L2}, we have
$$
\sigma_{1,4}(\b) \geq w_{1,4} = \alpha, \sigma_{2,4}(\b)\geq w_{2,4} = \alpha, \sigma_{1,3}(\b) \geq w_{1,3} = \beta, \sigma_{2,3}(\b) \geq w_{2,3} = \beta.
$$
It follows that
\begin{equation} \label{eq4}
\begin{cases}
\alpha \leq b_1 + b_3 \leq 2\alpha - \beta - 2 \\
1 \leq b_1 - b_3 \leq \alpha - \beta - 1.
\end{cases}
\end{equation}
In a similar manner to the Subcase 1.2, we find that the number of vectors $\b$ of the form \eqref{form1} with $0< b_1, b_3 < \alpha - 1$ and $b_1, b_3$ satisfying \eqref{eq4} equals
$$
\dfrac{1}{4}\biggl[2(\alpha - \beta - 1)^2 - \overline{\alpha - 2\beta + 1}^2\biggr].
$$
For such a vector $\b$, we get $\mathcal F(\Delta_\b(I))=\{\{1,2\}, \{3,4\}\}$ by Lemma~\ref{L2}, so 
$$
\dim_K H^{1}_{\mm}(R/I)_\b = \dim_K\widetilde{H}_{0}(\Delta_\b(I);K) = 1.
$$ 
Conversely, for each type-2 pair of edges $(ij, pq)$ of $G$, there are
$$
\dfrac{1}{4}\biggl[2(\alpha - \beta - 1)^2 - \overline{\alpha - 2\beta + 1}^2\biggr]
$$
vectors $\b$ of the form
$$
(0, \ldots, \underbrace{b_1}_{i}, \ldots, \underbrace{\alpha-1-b_1}_{j}, \ldots, \underbrace{b_3}_{p}, \ldots, \underbrace{\alpha - 1 - b_3}_{q}, 0, \ldots, 0)
$$
with $0 < b_1, b_3 < \alpha - 1$ and $b_1, b_3$ satisfying the system \eqref{eq4}. These vectors arise differently to those in the Subcases 1.1 and 1.2 by the configuration of $G[i, j, p, q]$. As a result, in the Subcase 1.3 the number of vectors $\b$ is 
\begin{equation} \label{case1.3}
\dfrac{1}{4}\biggl[2(\alpha - \beta - 1)^2 - \overline{\alpha - 2\beta + 1}^2\biggr] p_2(G).
\end{equation}
Again, for each of these vectors $\b$, we have $\dim_K H^{1}_{\mm}(R/I^{(m)})_{\b} = 1$.

From \eqref{case1.1}, \eqref{case1.2}, and \eqref{case1.3}, we conclude that in the Case 1, 
\begin{equation*} 
\begin{split}
\beta_{n-1, n + a_1}(R/I) = &\biggl[(\alpha - \beta)^2 + (\alpha - \beta - 1)^2 - \overline{\alpha - 2\beta + 1}^2\biggr]p_0(G) + \\
& \dfrac{1}{2}\biggl[(\alpha - \beta - 1)(2\alpha - 2\beta - 1) - \overline{\alpha - 2\beta + 1}^2\biggr] p_1(G) + \\
& \dfrac{1}{4}\biggl[2(\alpha - \beta - 1)^2 - \overline{\alpha - 2\beta + 1}^2\biggr] p_2(G) + \overline{\alpha - 2\beta}\, b(G).
\end{split}
\end{equation*}

\noindent {\bf Case 2. $G$ does not satisfy $(G_1)$}. In this case, $G$ must satisfy $(G_2)$ as $R/I$ is not Cohen-Macaulay. By \cite[Theorem 2.3]{MT}, $|\b| = \alpha + \beta - 2$. Moreover, the induced subgraph $G[1, 2, 3, 4]$ must have exactly one edge by \cite[Lemma 4.4]{MN}, say $\{1, 2\}$. Then
$$
|\b| = \alpha + \beta - 2 = w_{1,2} + w_{3,4} - 2 \geq \sigma_{1,2}(\b) + \sigma_{3,4}(\b) = |\b| + \sum_{i \geq 5} b_i \geq |\b|. 
$$
Thus $\sigma_{1,2}(\b) = \alpha - 1$, $\sigma_{3,4}(\b) = \beta - 1$ and $b_i = 0$ for all $i \geq 5$, so the vector $\b$ has the form 
\begin{equation} \label{form11}
(b_1, \beta - 1 - b_1, b_3, \alpha - 1 - b_3, 0, \ldots, 0).
\end{equation}
Moreover, we have
$$
\sigma_{1,3}(\b)\geq w_{1,3} = \beta, \sigma_{1,4}(\b) \geq w_{1,4} = \beta, \sigma_{2,3}(\b) \geq w_{2,3} = \beta, \text{ and } \sigma_{2,4}(\b) \geq w_{2,4} = \beta.
$$
It follows that
\begin{equation} \label{eq5}
\begin{cases}
\beta \leq b_1 + b_3 \leq \alpha - 2 \\
1 \leq b_3 - b_1 \leq \alpha - \beta - 1
\end{cases}.
\end{equation}

\noindent {\bf Subcase 2.1.} $\beta > 1$. Assuming $\alpha$ is odd and $\beta$ is even, then by Lemma~\ref{sol2}, the number of vectors $\b$ with $0 < b_1 < \beta - 1$, $0 < b_3 < \alpha - 1$ and $b_1, b_3$ satisfying \eqref{eq5} equals
$$
\dfrac{1}{2}\biggl[(\alpha - \beta - 1)^2 - \overline{\alpha - 2\beta + 1}^2\biggr].
$$
For such a vector $\b$, we have $\Delta_\b(I)$ consists of two connected components $\D_1, \D_2$ with $\D_1=\{\{1,2\}\}$ and $\{3,4\}\in\D_2$ by Lemma~\ref{L2}, so 
$$
\dim_K H^{1}_{\mm}(R/I)_\b = \dim_K\widetilde{H}_{0}(\Delta_\b(I);K) = 1.
$$ 
Conversely, for each $pq \not \in G$ such that $\{p, q\} \cup G$ is disconnected, there are
$$ 
\dfrac{1}{2}\biggl[(\alpha - \beta - 1)^2 - \overline{\alpha - 2\beta + 1}^2\biggr]
$$ 
vectors $\b$ of the form 
$$
(0, \ldots, \underbrace{b_1}_{i}, \ldots, \underbrace{\alpha-1-b_1}_{j}, \ldots, \underbrace{b_3}_{p}, \ldots, \underbrace{\alpha - 1 - b_3}_{q}, 0, \ldots, 0)
$$
with $0 < b_1 < \beta - 1, 0 < b_3 < \alpha - 1$ and $b_1, b_3$ satisfying the system \eqref{eq5}. 

If $b_1 = 0$ and $0 < b_3 < \alpha - 1$, then $\b = {\b'} = (0, \beta - 1, b_3, \alpha - 1 - b_3, 0, \ldots, 0)$. The system \eqref{eq5} turns into $\beta \leq b_3 \leq \alpha - \beta - 1$, so the number of vectors ${\b'}$ is $\overline{\alpha - 2\beta}$. Then $\mathcal F(\Delta_\b(I))=\{3,4\}\cup \{\{3,i\} ~|~3i\in G \}\cup \{\{4,i\} ~|~ 4i\in G \}\cup\{\{2,i\} ~|~ 2i\in G \}$, thus
$$
\dim_K H^{1}_{\mm}(R/I)_{\b'} = \dim_K\widetilde{H}_{0}(\Delta_{\b'}(I);K) = 1.
$$
The situation $b_1 = \beta - 1$ and $0 < b_3 < \alpha - 1$ is analogous. In general, for each $pq \not \in G$ such that $\{p, q\} \cup G$ is disconnected and for each $i \in V(G)$ (note that $G$ is connected, so $i$ belongs to some edge of $G$), there are $\overline{\alpha - 2\beta}$ vectors ${\b'}$ as above with $\beta-1$ in the $i$th coordinate and other two nonzero coordinates in the $p$th and $q$th.

If $b_3 \in \{0, \alpha-1\}$, then it follows from \eqref{eq5} that $b_1 \leq -1$, which is impossible. Consequently, the number of vectors $\b$ in Subcase 2.1 equals 
\begin{align*} 
\beta_{n-1, n + a_1}(R/I) &= \biggl[\dfrac{(\alpha - \beta - 1)^2}{2} - \dfrac{\overline{\alpha - 2\beta + 1}^2}{2} + \overline{\alpha - 2\beta} \, |V(G)|\biggr]\binom{n-|V(G)|}{2},
\end{align*}
as for each $\b$, we have $\dim_K H^{1}_{\mm}(R/I^{(m)})_{\b} = \dim_K\widetilde{H}_{0}(\Delta_{\b'}(I);K) = 1$. 

\medskip

\noindent {\bf Subcase 2.2.} $\beta  = 1$. Then $\b = (0, 0, b_3, \alpha - 1 - b_3, 0, \ldots, 0)$ with $1 \leq b_3 \leq \alpha - 2$. Thus $\mathcal F(\Delta_\b(I))=\{3,4\}\cup \{\{i,j\} ~|~ ij\in G \}$ by Lemma~\ref{L2}, so $\dim_K\widetilde{H}_{0}(\Delta_\b(I);K) \neq 0$ (i.e.,  $\dim_K\widetilde{H}_{0}(\Delta_\b(I);K) = 1$) if and only if $34\cup G$ is disconnected. Conversely, for each $pq \not \in G$ such that $\{p, q\} \cup G$ is disconnected, there are $\alpha - 2$ vectors $\b$ with the only two nonzero coordinates in the $p$th and $q$th.  Consequently, we get
$$
\beta_{n-1, n + a_1}(R/I) = (\alpha - 2)\binom{n-|V(G)|}{2}
$$
when $\beta = 1$.
\end{proof}

As a consequence of the above computations, we obtain the proof of Theorem~\ref{extremal2}.

\begin{proof}[Proof of Theorem~\ref{extremal2}]
We have pointed out that $R/I$ has two extremal Betti numbers if and only if it satisfies the condition $(G_1)$ or $(G_2)$ (i.e, $R/I$ is not Cohen-Macaulay) and $a_1 + 1 < a_2 + 2$. By \cite[Theorems 2.2 and 3.5]{MT} (where the values of $a_1$ and $a_2$ are computed), the latter happens in three cases: (i) $\girth(G) = 3$; (ii) $\girth(G) \neq 3$ and $G$ contains a vertex of degree $\geq 2$; (iii) $G$ consists of at least two disjoint edges and $\alpha < 2\beta$, or $G$ consists of exactly one edge. 

On the other hand, $\beta_{n-2, n + a_2}$ is the unique extremal Betti number of $R/I$ if and only if $R/I$ is Cohen-Macaulay, while $\beta_{n-1, n + a_1}$ is the unique extremal Betti number only when $a_1 + 1 \geq a_2 + 2$, that is, when $G$ consists of at least two disjoint edges and $\alpha \geq 2\beta$.	

In each case, the concrete values of $\beta_{n-2, n + a_2}$ are computed in Proposition~\ref{3}, while the values of $\beta_{n-1, n + a_1}$ are taken from Proposition~\ref{thm2}.
\end{proof}

\section{Applications}
In this section, we classify the rings $R/I$ which are pseudo-Gorenstein in terms of $(\alpha, \beta)$ and $G$, and give several examples illustrating our results.

The definition of pseudo-Gorenstein rings was introduced in \cite{EHHM} as a generali\-zation of Gorenstein property. Recall that $R/I$ is pseudo-Gorenstein if it is Cohen-Macaulay and $\beta_{n-2, n + a_2}(R/I) = 1$. In our case, this amounts to $\beta_{n-2, n + a_2}(R/I)$ is the unique extremal Betti number of $R/I$ and $\beta_{n-2, n + a_2}(R/I) = 1$.

As a direct consequence of Theorems~\ref{extremal1}--\ref{extremal5}, we obtain:

\begin{Corollary} $R/I$ is a pseudo-Gorenstein ring if and only if one of the following conditions holds:
\begin{enumerate}
\item[(i)] $(\alpha,\beta) = (2,1)$, $G$ has exactly one vertex of degree $2$ and other vertices of degree $1$. 
\item[(ii)] $(\alpha,\beta) = (2,1)$ and $G$ contains exactly one cycle of length $3$.
\item[(iii)] $(\alpha,\beta) = (4,3), s(G)=0$ and $G$ contains exactly one cycle of length $3$.
\item[(iv)] $\alpha=\beta+1,\beta \neq 1,3, \, p_0(G)=0$ and $G$ contains exactly one cycle of length $3$.
\item[(v)] $\alpha=\beta+1,\beta \neq 1$ and $G$ is the star $K_{1,2}$ .
\item[(vi)] $(\alpha,\beta) = (4,2), n=5$ and $G$ is either the cycle $C_3$ or the star $K_{1,2}$.
\end{enumerate}
\end{Corollary}

Next, we provide various examples to illustrate our results. The examples presented here have been inspired by computations performed by the computer algebra system \cite{M2}.
\begin{Example} Let $n=8$ and $G=\{12,13,23,24,34,45,56,67\}$. Then $G$ satisfies $(G_1)$ and $\girth(G)=3$. We observe that

- Disconnected pairs of edges are $$(12,56);(12,67);(13,56);(13,67);(23,56);(23,67);(24,67);(34,67).$$ So $p_0(G)=8$ and $G$ does not satisfy $(G_3)$. 

- Type-1 pairs of edges are $(12,45);(13,45);(24,56);(34,56);(45,67)$. So $p_1(G)=5$.

- There is only a type-2 pair of edges which is $(23,45)$.  So $p_2(G)=1$.

\noindent Moreover $b(G)={\underset{pq\in G}\sum} \, |V(G) \backslash (N(p)\cup N(q))|=22$ and $c_3(G)=2$. 

Let $I\in \C_8(\alpha,\beta)$ and $w_{i,j}=\alpha$ if and only if $ij\in G$. By Theorem \ref{extremal1}, Theorem \ref{extremal2}, Theorem \ref{extremal3}, Theorem \ref{extremal4} and Theorem \ref{extremal5}, $R/I$ always has two extremal Betti numbers. Moreover, the first extremal Betti number $\beta_{6,a_2}(R/I)$ is always equal to $c_3(G)=2$ for any $\alpha>\beta>0$. The other extremal Betti number is determined as follows.
\begin{enumerate}
\item Let $\alpha=8$. Then, the second extremal Betti number is
$$
\begin{cases}
\beta_{7,22}(R/I) =B_1(8,3) = 376&\text{ if } \beta=3 \text { by Theorem \ref{extremal2}},\\
 \beta_{7,21}(R/I) =(n-4)p_0(G)= 32 &\text{ if } \beta=7 \text { by Theorem \ref{extremal1}},\\
 \beta_{7,21}(R/I) =4p_0(G)+p_1(G)= 37&\text{ if } \beta=6 \text { by Theorem \ref{extremal4}}.\\
\end{cases}
$$
\item Let $\alpha=7, \beta =5$. Then, by Theorem \ref{extremal3}, the second extremal Betti number is $$\beta_{7,20}(R/I) =5p_0(G)+2p_1(G)+p_2(G) = 51.$$ 
\item Let $\alpha=4, \beta =2$. Then, by Theorem \ref{extremal5}, the second extremal Betti number is $$\beta_{7,14}(R/I) =4p_0(G)+p_1(G) = 37.$$

\end{enumerate}
\end{Example}

Finally, we end this article with a conjectural observation. This also suggests that studying extremal Betti numbers helps in understanding Betti tables of the defining ideal of hyperplanes $L_{i,j}$ with respect to the weights $w_{i,j}$.

\begin{Remark}
Let $I\in \C_n(\alpha,\beta)$ for $n \geq 5$. When $\beta_{n-2, n+a_2}(R/I)$ is an extremal Betti number and $\beta \geq 2$, our computations using Macaulay2 suggest that the last row of the Betti diagram of $R/I$ would be symmetric; moreover, it should stabilize when the graph $G$ is fixed and $\alpha > \beta \geq 2$ vary.  Note that for $\beta \geq 2$, $\beta_{n-2, n+a_2}(R/I)$ is not extremal only when $G$ consists of at least two disjoint edges, and $\alpha \geq \max\{\beta + 3, 2\beta\}$ or $(\alpha, \beta) = (4,2)$.
\end{Remark}

\subsection*{Acknowledgments} Part of this work was done while the last two authors were visiting Vietnam Institute for Advanced Study in Mathematics (VIASM).  We would like to thank VIASM for its hospitality and financial support. 

\end{document}